\documentclass[11pt,a4paper]{article}

\usepackage[english]{babel}
\usepackage[utf8]{inputenc}
\usepackage{amssymb}
\usepackage{amsfonts}
\usepackage{amsmath}
\usepackage{amsthm}
\usepackage{graphicx}
\usepackage[colorinlistoftodos]{todonotes}
\usepackage{authblk}
\usepackage{cite} 
\usepackage[labelfont=bf,labelsep=space]{caption}
\usepackage[framed,numbered,autolinebreaks,useliterate]{mcode}

\usepackage{CJK}
\usepackage{color}
\usepackage{epsfig}
\usepackage{mathptmx}
\usepackage{mathrsfs}
\usepackage{booktabs}
\usepackage{graphicx}
\usepackage{epstopdf}
\usepackage{multirow}
\usepackage[subfigure]{graphfig}
\usepackage{graphics}
\usepackage{times}
\usepackage{booktabs}
\usepackage{hyperref}
\usepackage{fancyhdr}

\newtheorem{theorem}{Theorem}[section]

\newtheorem{coro}{Corollary}[section]

\newtheorem{remark}{Remark}[section]

\begin{document}

\title{Design and Validation of a Virtual Player for Studying Interpersonal Coordination in the Mirror Game}

\author[1]{Chao Zhai}
\author[1]{Francesco Alderisio}
\author[2]{Piotr S\l{}owi\'{n}ski}
\author[2]{Krasimira Tsaneva-Atanasova}
\author[1,3]{Mario di~Bernardo}

\affil[1]{Department of Engineering Mathematics, University of Bristol, Merchant Venturers' Building, BS8 1UB, United Kingdom}
\affil[2]{College of Engineering, Mathematics and Physical Sciences, University of Exeter, United Kingdom}
\affil[3]{Department of Electrical Engineering and Information Technology, University of Naples Federico II, 80125 Naples, Italy}
\maketitle

\begin{abstract}
The mirror game has been recently proposed as a simple, yet powerful paradigm for studying interpersonal interactions. It has been suggested that a virtual partner able to play the game with human subjects can be an effective tool to affect the underlying neural processes needed to establish the necessary connections between the players, and also to provide new clinical interventions for the rehabilitation of patients suffering from social disorders. Inspired by the motor processes of the central nervous system (CNS) and the musculoskeletal system in the human body, in this paper we develop a novel interactive cognitive architecture based on nonlinear control theory to drive a virtual player (VP) to play the mirror game with a human player (HP) in different configurations. Specifically, we consider two cases: the former where the VP acts as leader and the latter where it acts as follower. The crucial problem is to design a feedback control architecture capable of imitating and following or leading a human player (HP) in a joint action task. Movement of the end-effector of the VP is modeled by means of a feedback controlled Haken-Kelso-Bunz (HKB) oscillator, which is coupled with the observed motion of the HP measured in real time. To this aim, two types of control algorithms (adaptive control and optimal control) are used and implemented on the HKB model so that the VP can generate  human-like motion while satisfying certain kinematic constraints. A proof of convergence of the control algorithms is presented in the paper together with an extensive numerical and experimental validation of their effectiveness. A comparison with other existing designs is also discussed, showing the flexibility and the advantages of our control-based approach.
\end{abstract}

\section{Introduction}
The emergence of coordinated behavior between humans is a common phenomenon in many areas of human endeavor. Examples include  improvisation theater, group dance, music playing, team sports and parade marching \cite{john87}. At the core of the interaction between the players lies a fundamental feedback mechanism where each player adapts his/her motion in response to the observed movement of the other.

To study this intriguing phenomenon, the mirror game has been recently proposed as a simple, yet effective paradigm. In its simplest formulation, the mirror game features two people imitating each other's movements at high temporal and spatial resolution \cite{noy11}. The game can be played in different experimental conditions: the former where one of the players leads and the other has to follow the leader movement (Leader-Follower condition); the latter where the two players create joint synchronized movement (Joint Improvisation condition).

The theory of similarity in social psychology suggests that people prefer to team up with others possessing similar morphological and behavioral features, and that interpersonal coordination is enhanced if their movement shares similar kinematic features  \cite{fol82, lak11}. Further evidence suggests that motor processes caused by interpersonal coordination are strictly related to mental connectedness. To be specific, motor coordination between two people contributes to social attachment \cite{wil09}.

As suggested in \cite{wil09}, coordination games can therefore be used to help people suffering from social disorders to improve their social skills. Also they can be effectively exploited in social robotics to enhance attachment, coordination and rehabilitation during human-robot interactions \cite{hussain}.
For this reason, it has  been proposed that creating a VP or avatar able to coordinate its motion with that of a HP can be extremely useful to study the onset of coordination and how it is affected by similarity/dissimilarity between the players' motion characteristics \cite{pnas14}. A VP can also be used for diagnostics and rehabilitation of patients suffering from social disorders as recently proposed in \cite{alter}.

The aim of this paper is the design of a novel interactive cognitive architecture (ICA) based on nonlinear control theory able to drive a VP to play the mirror game with a human either as a leader or as a follower. Specifically, the goal is that of designing a cognitive architecture able to drive the motion of the VP interacting with a HP in real-time while exhibiting certain desired kinematic features. When playing as a follower, the ICA needs to guarantee that, while exhibiting the desired movement properties, the VP tracks as closely as possible the motion of the human leader. When playing as the leader, the ICA needs instead to generate new interesting motion. In both cases, it is crucial for the VP to engage with the HP by producing human-like response in terms of kinematics (maximum acceleration, velocity profile etc) and delay times. In this paper we take the view that the design of such an architecture is fundamentally a nonlinear control design problem where given some reference input the architecture has to drive the VP onto a desired motion  which is a function of the movement of the human player being sensed during the game. In particular, the ICA can be integrated into the humanoid robot to achieve the desired dual-arm coordination \cite{liu}.

We explore two different approaches, one based on adaptive control, the other on optimal control. Our control architecture mimics the two fundamental actions which have been suggested to be at the core of the emergence of  motor coordination between two or more effectors in biological systems: feedback and feedforward \cite{jor99,drop,pool}. Specifically, the motor system is able to correct the deviation from the desired movement with the aid of feedback control, whilst feedforward control allows it to reconcile the interdependency of the involved effectors and preplan the response to the sensory incoming information \cite{drop,pool,des00}.

It is shown experimentally that the proposed control architectures are able to effectively drive the VP to play the mirror game while generating motion with desired kinematic properties. In particular, we use the concept of Individual Motor Signature (IMS) recently proposed in \cite{piotr,piotr15} to characterize the motion of an individual player and evaluate how similar/dissimilar the motion of two different individuals is. Following our approach we are able to show that the VP driven by the cognitive architecture presented in the rest of this paper can play the mirror game either as a leader or a follower while exhibiting a desired IMS.

Relevant previous work in the literature includes the generation of human-like movement \cite{zhang14}, the development of a mathematical model to explain the coordination dynamics observed experimentally in the mirror game \cite{noy11}, and the Human Dynamic Clamp paradigm proposed in \cite{pnas14,kel09,hol09} where the use of a virtual partner driven by appropriate mathematical models is proposed to study human motor coordination. These previous approaches will be used to investigate and compare the performance of the novel strategy presented in this paper. We wish to emphasize that the control algorithms developed and validated in what follows can be also effectively used for trajectory planning to enhance human-robot coordination in joint interactive tasks.

The rest of the paper is organized as follows. The mirror game set-up, problem statement and motor signature are discussed in Section \ref{sec:problem} before presenting the schematic of the proposed cognitive architecture in Section \ref{sec:ca}. The feedback control strategies at the core of the ICA are developed and analyzed in Section \ref{sec:adaptive} and \ref{sec:optimal}. The experimental validation of the control algorithms is presented in Section \ref{sec:validation} where experimental results are discussed showing the effectiveness of the proposed strategies. A comparison with other existing approaches is also carried out.  Finally, conclusions and suggestions for future work are drawn in Section \ref{sec:conclusions}.

\section{The Mirror Game Problem}\label{sec:problem}

Investigation of interpersonal coordination requires appropriate experimental paradigms. A typical paradigm recently proposed in the literature is the mirror game, which involves two people imitating each other's movements at high temporal and spatial resolution \cite{noy11}. It can be played in two different conditions: \emph{Leader-Follower condition}, where the follower attempts at tracking the leader motion as accurately as possible, and \emph{Joint Improvisation condition}, where the players jointly coordinate and synchronize their movements without any of the two being designated as leader or follower.

\begin{figure}
\scalebox{0.25}[0.25]{\includegraphics{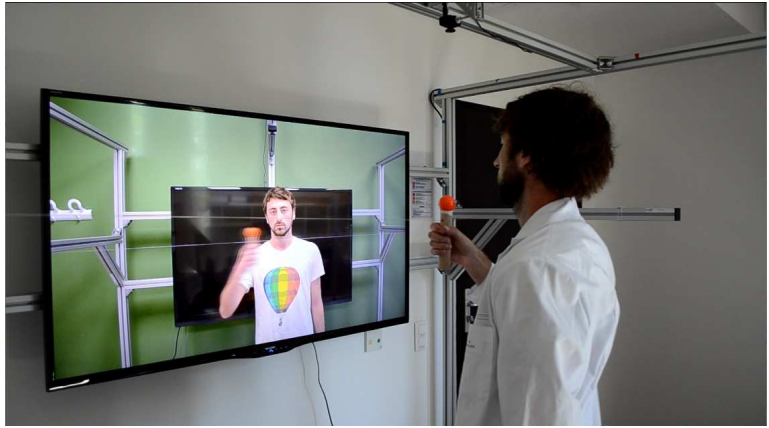}}\centering
\caption{\label{setup}Experimental set-up of the mirror game between a VP and a HP at the University of Montpellier, France (see \cite{alter} for further details).}
\centering
\end{figure}

Our set-up is inspired by the one in \cite{noy11}. Specifically, a small orange ball is mounted onto a string, which the HP can move back and forth along the string itself. In the meanwhile, the VP on the opposite screen moves its own ball on a parallel string with the same length (see Fig. \ref{setup}). In this implementation of the mirror game, two players (a HP and a VP) are required to move their respective ball back and forth and synchronize their movement. Here, we assume that the game is played in a Leader-Follower condition, where the HP is the leader and the VP (robot or computer avatar) is the follower trying to track the leader movement. However, the VP can opt to act as the leader as well.

The position of the ball moved by the HP is detected by a camera. A feedback control strategy then needs to be designed in order to generate the trajectory of the ball moved by the VP so as to track the movement of the ball controlled by the HP. Such a trajectory can then be provided to the on-board controllers of the VP  (robot or computer avatar) as the desired trajectory for its end effector.

To solve this control problem so that the VP motion presents similar features to the motion of a human player, we need to choose an appropriate model of the VP motion that can then be controlled using a nonlinear feedback strategy. To this purpose, here we use the Haken-Kelso-Bunz oscillator, which was first proposed in \cite{hkb85} as a model able to capture the observations made in experiments on human bimanual coordination. The model consists of two nonlinearly coupled nonlinear oscillators described by
\begin{equation}
\label{eq:HKB}
\ddot{z}+(\alpha \dot{z}^2+\beta z^2-\gamma)\dot{z}+\omega^2z=[a+b(z-w)^2](\dot{z}-\dot{w})
\end{equation}
where $z, \dot{z}$ represent the position and velocity of finger $1$, $w, \dot w$ the position and velocity of finger $2$ (modeled by a replica of the equation above obtained by swapping $w, \dot w$ with $z, \dot z$); $a$ and $b$ are the coupling parameters and $\alpha$, $\beta$, $\gamma$ and $\omega$ characterize the response of each uncoupled finger when subject to some reference signal. However, it is worth pointing out that, other than describing intrapersonal motor coordination, the HKB model has been also used to describe interpersonal motor coordination involving two different people \cite{fuchs96,schmidt08}. In particular, the HKB oscillator has been suggested in the literature as a paradigmatic example of human motor coordination \cite{pnas14,hkb85}. Solving the mirror game can then be formulated as the following control problem. Given a nonlinear HKB oscillator of the form
\begin{equation}\label{system}
\left\{
  \begin{array}{ll}
    \dot{x}=y \\
    \dot{y}=-(\alpha y^2+\beta x^2-\gamma)y-\omega^2x+u
  \end{array}
\right.
\end{equation}
where $x$ and $\dot{x}$ refer to the position and velocity of the end effector of the VP, respectively, and $u$ is an external control input, the problem is to design a feedback controller $u$ such that $x(t)$ achieves bounded asymptotic tracking of the position of the HP, while expressing some desired kinematic features.

\begin{figure}\centering
 \subfigure[Same player]
 {\includegraphics[width=0.7\textwidth]{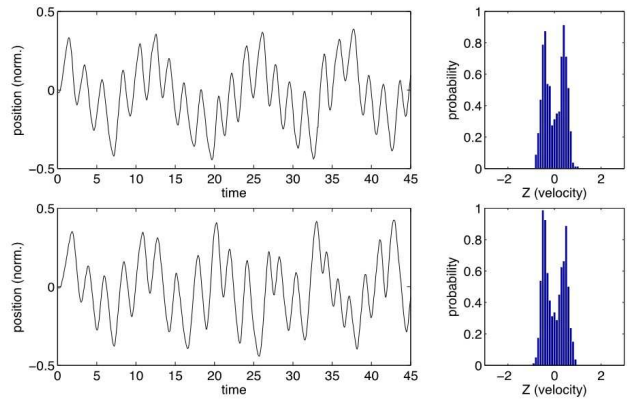}}
 \subfigure[Two different players]
 {\includegraphics[width=0.7\textwidth]{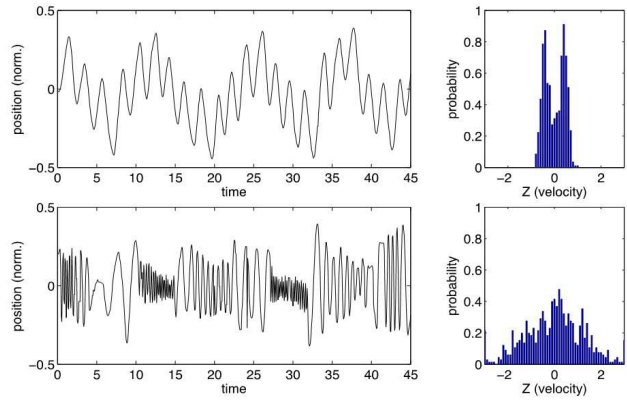}}
 \caption{Position time series and the corresponding PDF of velocity time series for the same player (a) and two different players (b).}
 \label{pdf_vel}
\end{figure}

As metrics to characterize the kinematic properties of the motion of an individual playing the game  we use the concept of individual motor signature (IMS), recently introduced in \cite{piotr,piotr15}. It has been shown that the IMS is time invariant and is unique for each player. It is defined in terms of the velocity profile (or distribution) of the player's motion in solo trials. To quantify how similar or dissimilar the signatures of two different players are, we use the earth mover's distance (EMD) between any two probability distribution functions (PDF) of their velocity time series \cite{piotr15,levi01}. The EMD between two PDFs $p_1$ and $p_2$ can be computed as follows
$$
EMD(p_1,p_2)=\int_{Z}|CDF_{p_1}(z)-CDF_{p_2}(z)|dz
$$
where $Z$ denotes the integration domain, and $CDF_{p_i}(z)$ denotes the cumulative distribution function of the distribution $p_i,i\in\{1,2\}$. Fig.\ref{pdf_vel}(a) shows the position time series of the same HP and the corresponding PDF of velocity in two solo trials. It is visible that the two PDFs of velocity time series resemble each other in terms of their shape, and the EMD between them is $0.024$. In contrast, the two PDFs of velocity time series in Fig.\ref{pdf_vel}(b) differ remarkably from each other, and the value of EMD is $0.604$, which confirms the qualification of the PDF of velocity time series in solo trials as individual motor signature.

\section{Design of the Cognitive Architecture}\label{sec:ca}

\begin{figure}
\scalebox{0.07}[0.07]{\includegraphics{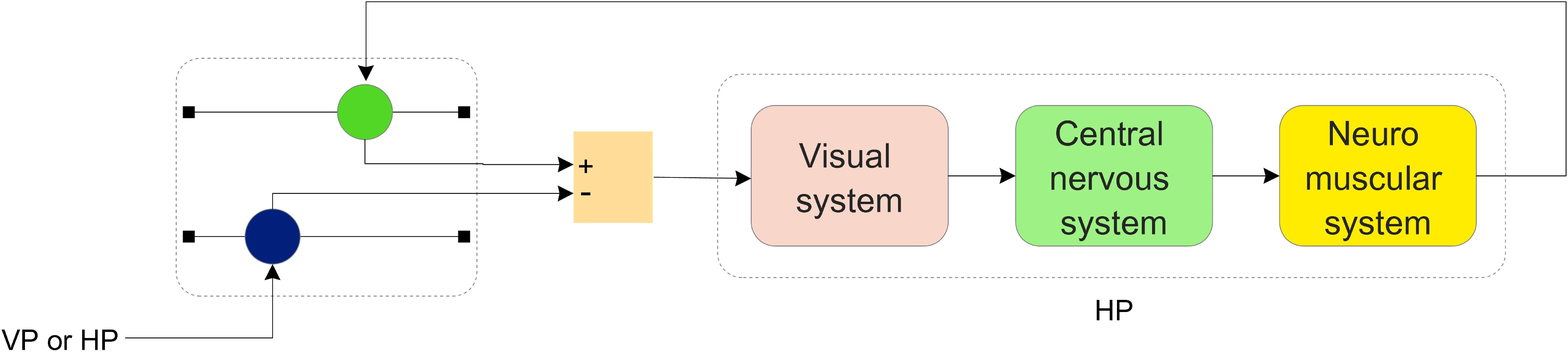}}\centering
\caption{\label{human}Motor coordination between two players in mirror game.}
\end{figure}

We design the cognitive architecture of the VP so as to replicate the main processes involved in making a human being play the mirror game (see Fig. \ref{human}). The visual system detects the ball's position on the string and generates visual signals, which are then transmitted to the central nervous system (CNS including brain and spinal cord). Several parts of the CNS (such as ventral horn, cerebellum and motor cortex) use an internal model to predict the kinematic characteristics of the other player's motion and generate the neural impulses that control the extension and contraction of muscles. Finally, the neuromuscular system activates and coordinates the muscles involved in generating the hand movements.

\begin{figure*}
\scalebox{0.055}[0.055]{\includegraphics{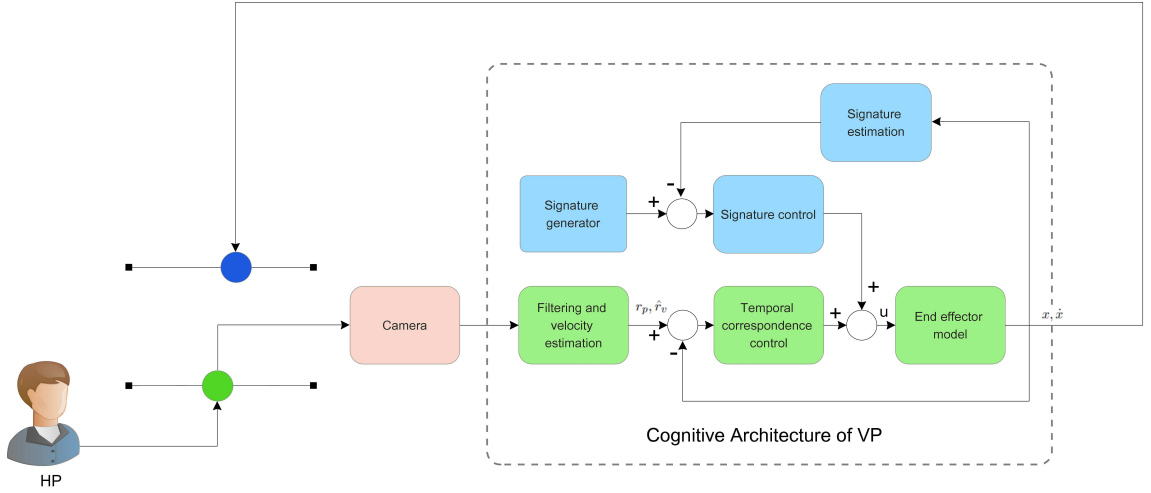}}\centering
\caption{\label{sig} Block diagram of the cognitive architecture of VP. Green blocks allow for the control of temporal correspondence between the VP and the HP, and blue blocks take into account the desired individual motor signature.}
\end{figure*}

This architecture is mapped onto the real-time control schematics shown in Fig. \ref{sig} whose blocks are briefly described below.

\begin{itemize}
  \item A {\em camera} is used to detect the position of the HP, say $r_p$;
  \item A {\em filtering and velocity estimation block} is used to filter the position data acquired by the camera via a low pass filter and to estimate the velocity of the HP (reference) via the simple formula
\begin{equation}\label{rvest}
\hat{r}_v(t)=\frac{r_p(t_k)-r_p(t_{k-1})}{T} \quad t\in[t_k, t_{k+1}]
\end{equation}
where $k\in N^*$, and $T=t_k-t_{k-1}$ denotes the sampling period of the camera.
The estimated velocity is then used to predict the HP position over the next interval by using the expression:
\begin{equation}\label{rpest}
\hat{r}_p(t)=r_p(t_k)+\hat{r}_v(t)(t-t_k), \quad t\in[t_k, t_{k+1}]
\end{equation}

 As an alternative, we could adopt a nonlinear observer to provide a better prediction of the reference velocity; for example, the nonlinear extended observer in \cite{wan03}. Here we find that such a complication is unnecessary to solve the problem of interest and therefore choose to use the simple yet effective estimation strategy discussed above.
  \item At the core of the architecture lie the two blocks {\em Temporal Correspondence Control}  and {\em Signature Control}. The former  is designed to regulate the end effector model so that its motion tracks that of  the HP with varying degrees of dynamic similarity. Specifically, it aims at minimizing the position error between the time series of the HP and that of the VP. The latter block uses the prerecorded velocity time series of a reference HP with the desired IMS (velocity profile) in order to generate the avatar trajectory with desired kinematic features. In particular, the aim of the signature controller is that of  reducing the distance (computed in terms of EMD) between the velocity distribution of the VP and that of some reference HP it aims at replicating the motion characteristics of.
  \item The prerecorded velocity trajectory of a reference HP playing solo representing the desired IMS is stored in the {\em Signature generator block} while the signature of the avatar motion is estimated by the {\em Signature estimation block}.
  \item The {\em end effector model} is used to generate the avatar motion via an appropriate feedback control scheme. As mentioned before, we use the HKB oscillator to describe the dynamics of the end effector.
  \item Finally, the output of the cognitive architecture (position and velocity $x$ and $\dot{x}$) is used as the reference motion for the VP.
\end{itemize}

In what follows we focus on the design of the feedback control strategies that drive the cognitive architecture. We derive and compare two different types of controllers. First we develop an adaptive algorithm able to control the temporal correspondence between the VP and the HP during the game (green blocks in Fig.\ref{sig}). Then, we consider an optimal controller to solve simultaneously the multi objective control problem of tracking the trajectory of the HP while preserving the features of the desired IMS of interest (both green and blue blocks in Fig.\ref{sig}).
For both strategies a proof of convergence is given before presenting numerical  and experimental investigation of their performance.

\section{Adaptive Control of Temporal Correspondence} \label{sec:adaptive}

To solve the control problem of temporal correspondence, we propose an adaptive controller based on the end effector model shown in (\ref{system}). Specifically, we choose the nonlinear controller given by
\begin{equation}
\label{eq:u}
u=\underbrace{[a(t)+b(t)(x-r_p)^2](\dot{x}-\hat{r}_v)}_{Coordination}-\underbrace{C_pe^{-\delta (\dot{x}-\hat{r}_v)^2}(x-r_p)}_{Temporal~Correspondence}
\end{equation}
where $r_p$ is the position of the HP, $\hat r_v$ is the estimated velocity, $C_p$ and $\delta$ are constant parameters while the coupling parameters $a$ and $b$ are updated according to the adaptive laws:
\begin{equation}\label{coup:a}
\begin{split}
\dot{a}&=-e^{-2a} \left[ (x-r_p)(y-\hat{r}_v)+\eta_a(x-r_p)^2 \right]-\eta_a
\end{split}
\end{equation}
and
\begin{equation}\label{coup:b}
\begin{split}
\dot{b}&=  \frac{y-\hat{r}_v}{e^{2b}}[\omega^2x+(\alpha y^2+\beta x^2-\gamma)y-\eta_a(y-\hat{r}_v)-u]-\eta_a
\end{split}
\end{equation}
where $\eta_a$ is a positive constant. Note that the control law (\ref{eq:u}) consists of two complementary terms. The first has the same structure as the one of the coupling proposed in \cite{hkb85} to model the interaction between two HPs, albeit with the introduction of adaptive parameters to account for variability between different HPs. The second term, depending on the fixed parameters $C_p$ and $\delta$, deals with the position error when the velocity mismatch approaches zero and hence the first term decays to zero. When $|\dot{x}-\hat{r}_v|$ is relatively large, the coupling term of the HKB equation instead dominates and motor coordination between the two players becomes more pronounced during the mirror game.

Theoretical analysis of the adaptive control algorithm in Table \ref{table_ac} is given in what follows below.

\subsection{Convergence analysis}

\begin{table}
 \caption{\label{table_ac} Adaptive Control Algorithm}
 \begin{center}
 \begin{tabular}{lcl}
  \midrule
  1: set $k=1$ and running time $T_s$ \\
  2: \textbf{while} ($time<T_s$) \\
  3: ~~~~~~~detect the position $r_p(t_k)$ of HP  \\
  4: ~~~~~~~estimate the velocity $\hat{r}_v(t_k)$ of HP  with (\ref{rvest}) \\
  5: ~~~~~~~generate the control signal $u$ with (\ref{eq:u}) \\
  6: ~~~~~~~update coupling parameters $a$ and $b$ with (\ref{coup:a}) and (\ref{coup:b})  \\
  7: ~~~~~~~obtain the position $x$ and velocity $\dot{x}$ of VP by solving (\ref{system})\\
  8: ~~~~~~~$k=k+1$\\
  9: \textbf{end while} \\
  \bottomrule
 \end{tabular}
 \end{center}
\end{table}

\begin{theorem}
The adaptive feedback controller (\ref{eq:u}) ensures that the solution of the controlled HKB model (\ref{system}) satisfies
$$
|x(t)-r_p(t)|\leq e^{\eta_a T}\sqrt{\frac{2\epsilon}{e^{2\eta_a T}-2}}+\frac{2}{e^{\eta_a T}}\sqrt{E(0)}, \quad t\in [0,+\infty)
$$
if $\eta_a$ is chosen so that
$$
\eta_a>\frac{\ln 2}{2T}
$$
where
$$
E(0)=\frac{1}{2} \left[ (x(0)-r_p(0))^2+(y(0)-\hat{r}_v(0))^2+e^{2a(0)}+e^{2b(0)} \right]
$$
and
$$
\epsilon=\sup_{k\in \mathbb{N}^*} (T^2+1)(\hat{r}_v(kT)-\hat{r}_v(kT-T))^2.
$$
\end{theorem}
Proof: Choose the energy-like function
\begin{equation}\label{eq:h}
E\triangleq\frac{1}{2} \left[ (x-r_p)^2+(y-\hat{r}_v)^2+e^{2a}+e^{2b}\right]
\end{equation}
Note that $\hat{r}_v$ is fixed in each sampling interval $[kT, (k+1)T), k\in \mathbb{N}^*$. Then the time derivative of $E$ along the trajectories of  (\ref{system}) with $u$ defined in (\ref{eq:u}) is given by
\begin{equation*}
\begin{split}
\dot{E}&=(x-r_p)(\dot{x}-\hat{r}_v)+(y-\hat{r}_v)\dot{y}+e^{2a}\dot{a}+e^{2b}\dot{b}\\
&=(x-r_p)(y-\hat{r}_v)-(y-\hat{r}_v)\left[(\alpha y^2+\beta x^2-\gamma)y+\omega^2x-u\right]\\
&+ \left( y-\hat{r}_v \right) [\omega^2x+(\alpha y^2+\beta x^2-\gamma)y-\eta_a(y-\hat{r}_v)-u]-\eta_a e^{2b}\\
&- \left[ (x-r_p)(y-\hat{r}_v)+\eta_a(x-r_p)^2 \right]-\eta_a e^{2a}\\
&=-\eta_a(x-r_p)^2-\eta_a(y-\hat{r}_v)^2-\eta_a e^{2a}-\eta_a e^{2b}\\
&=-2\eta_a E, \quad t\in [kT, (k+1)T)
\end{split}
\end{equation*}
Solving the above differential equation yields
\begin{equation}\label{eht}
E(t)=e^{-2\eta_a(t-kT)}E(kT),\quad t\in [kT,(k+1)T)
\end{equation}
Moreover, at the sampling point $kT$ we have
\begin{equation*}
\begin{split}
&~~~~E(kT)-E^-(kT)\\
&=\frac{1}{2}[(x-r_p(kT))^2-(x-r_p(kT-T)-\hat{r}_v(kT-T)T)^2\\
&+(y-\hat{r}_v(kT))^2-(y-\hat{r}_v(kT-T))^2]\\
&\leq (r_p(kT)-r_p(kT-T)-\hat{r}_v(kT-T)T)^2\\
&+(\hat{r}_v(kT)-\hat{r}_v(kT-T))^2+E^-(kT)\\
%&=(\hat{r}_v(kT)-\hat{r}_v(kT-T))^2T^2+(\hat{r}_v(kT)-\hat{r}_v(kT-T))^2+E^-(kT)\\
&=(1+T^2)(\hat{r}_v(kT)-\hat{r}_v(kT-T))^2+E^-(kT)\\
\end{split}
\end{equation*}
which is equivalent to
\begin{equation}\label{ineq}
E(kT)\leq \epsilon+2E^-(kT)
\end{equation}
where
$$
\epsilon=\sup_{k\in \mathbb{N}^*} (T^2+1)(\hat{r}_v(kT)-\hat{r}_v(kT-T))^2
$$
and
$$
E^-(kT)=\lim_{t\searrow kT}E(t)
$$
Evaluating (\ref{eht}) and (\ref{ineq}) at $t=(k+1)T$ and nesting the inequalities backwards till $t=0$, we get
\begin{equation}\label{inh}
\begin{split}
E(kT)&\leq\epsilon \left[ 1+\frac{2}{e^{2\eta_a T}}+ \left( \frac{2}{e^{2\eta_a T}} \right)^2+...+\left( \frac{2}{e^{2\eta_a T}} \right)^{k-1} \right] \\
&+ \left( \frac{2}{e^{2\eta_a T}} \right)^kE(0)\\
&=\epsilon\frac{1-(\frac{2}{e^{2\eta_a T}})^k}{1-\frac{2}{e^{2\eta_a T}}}+\left( \frac{2}{e^{2\eta_a T}} \right)^kE(0)\\
&\leq \frac{\epsilon}{1-\frac{2}{e^{2\eta_a T}}}+\frac{2}{e^{2\eta_a T}}E(0), \quad \forall k\in \mathbb{N}^*
\end{split}
\end{equation}
when the inequality $\eta_a >\frac{\ln 2}{2T}$ holds.
Moreover, combining (\ref{eht}) with (\ref{inh}), we get
$$
E(t)\leq\frac{\epsilon}{1-\frac{2}{e^{2\eta_a T}}}+\frac{2}{e^{2\eta_a T}}E(0),\quad t\in[0,+\infty)
$$
which clearly implies
$$
|x(t)-r_p(t)|\leq e^{\eta_a T}\sqrt{\frac{2\epsilon}{e^{2\eta_a T}-2}}+\frac{2}{e^{\eta_a T}}\sqrt{E(0)}, \quad t\in [0,+\infty)
$$
$\hfill{\blacksquare}$

\begin{remark}
It is easy to demonstrate that the coupling parameters $a$ and $b$ are upper bounded with the proposed adaptive laws.
\end{remark}

\begin{remark}
Since $r_p(t)\in[0,l], \forall t\geq0$ and $|\hat{r}_v(t)|\leq\frac{l}{T}$, the following inequality holds
\begin{equation*}
\begin{split}
\epsilon=\sup_{k\in \mathbb{N}^*} (T^2+1)(\hat{r}_v(kT)-\hat{r}_v(kT-T))^2\leq \frac{4l^2(1+T^2)}{T^2}
\end{split}
\end{equation*}
where $l$  refers to the length of the string. Generally, the upper bound for the position error is relatively conservative. When the velocity of the HP is small, $\epsilon$ is small as well, and the estimation for the position error is accurate enough.
In addition, taking the limit of \eqref{inh} as $kT \rightarrow \infty$ and combining it with \eqref{eq:h}, the position error between the two players satisfies the following inequality:
$$
\limsup_{kT \rightarrow +\infty}|x(t)-r_p(t)|\leq e^{\eta_a T}\sqrt{\frac{2\epsilon}{e^{2\eta_a T}-2}} %\quad t\rightarrow +\infty
$$
Similarly, we can estimate the velocity error as
$$
|\dot{x}(t)-\hat{r}_v(t)|\leq e^{\eta_a T}\sqrt{\frac{2\epsilon}{e^{2\eta_a T}-2}}+\frac{2}{e^{\eta_a T}}\sqrt{E(0)}, \quad t\in [0,+\infty)
$$
\end{remark}

While being effective in achieving bounded tracking of the HP position, the control approach derived so far is unable to explicitly guarantee that the generated motion follows a desired IMS. Therefore we consider a different scheme based on optimal control.
%%%%%%%%%%%%%%%%%%%%%%%%%%%%%%%%%%%%%%%%%%%%%%%%%%%%%%

\section{Optimal Temporal Correspondence and Signature Control}\label{sec:optimal}

\begin{figure}
\scalebox{0.06}[0.06]{\includegraphics{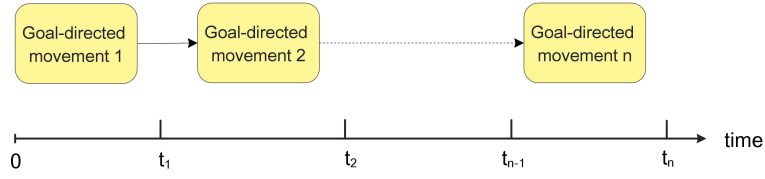}}\centering
\caption{\label{goalmov}Movement of the VP end effector in the mirror game}
\centering
\end{figure}

The second approach we propose is based on optimal control, which is an effective framework to allow for movement coordination and reconcile target tracking and individual motor signature \cite{died10}. We assume that in the mirror game the motion of the VP can be decomposed into a series of goal-directed movements (see Fig. \ref{goalmov}) influenced by both the position of the HP and the desired individual motor signature. Thus, we formulate the problem of driving the end effector motion as described by (\ref{system}) on a finite time interval $[t_k,t_{k+1}]$ as the dynamic optimization problem
\begin{equation}\label{minj}
\min_{u\in R} J
\end{equation}
where
\begin{equation}\label{vcost}
\begin{split}
J&=\frac{\theta_p}{2}\underbrace{(x(t_{k+1})-\hat{r}_p(t_{k+1}))^2}_{Temporal~Correspondence}\\ &+\frac{1}{2}\int_{t_k}^{t_{k+1}}\theta_{\sigma}\underbrace{(\dot{x}(\tau)-r_{\sigma}(\tau))^2}_{Similarity}+\eta_m u(\tau)^2d\tau
\end{split}
\end{equation}
with the constraint $\theta_p+\theta_{\sigma}=1$ and $\theta_p,\theta_{\sigma},\eta_m>0$ being tunable control parameters. Here, $\hat{r}_p(t_{k+1})$ denotes the estimated position of the HP at time $t_{k+1}$ (see (\ref{rpest}) for further details), while $r_{\sigma}$ refers to a prerecorded velocity time series representing the desired motor signature. For the sake of simplicity, the optimization interval $[t_k,t_{k+1}]$ in the cost function (\ref{vcost}) corresponds to the sampling interval. The above cost function mainly consists of three terms. The first term aims at minimizing the mismatch between the position time series of the HP and that of the VP. The second term takes care of making the velocity profile of the VP as close as possible to the reference one (motor signature). The last term guarantees boundedness of the control effort. In particular, the idea behind this cost function is that the human-like movement of the VP emerges from the integration of three different goals related to temporal correspondence, motor signature and control energy expenditure, respectively. Notice that the VP acts as a leader when $\theta_p$ is close to $0$, since the term related to the position error $x-\hat{r}_p$ in the cost function is negligible and the only aim of the VP is to exhibit the  motion characterized by the desired motor signature. Instead, the VP behaves as a follower if $\theta_p$ is close to $1$ as in this case the controller aims solely at minimizing the mismatch between the HP and VP terminal positions.

\subsection{Convergence Analysis}

\begin{table}
 \caption{\label{table_mpc} Optimal Control Algorithm}
 \begin{center}
 \begin{tabular}{lcl}
  \midrule
  1: set $k=1$ and running time $T_s$ \\
  2: \textbf{while} ($time<T_s$) \\
  3: ~~~~~~~detect the position $r_p(t_k)$ of HP  \\
  4: ~~~~~~~estimate the position $\hat{r}_p(t_{k+1})$ of HP with (\ref{rpest}) \\
  5: ~~~~~~~generate the control signal $u$ by solving (\ref{minj}) \\
  6: ~~~~~~~obtain the position $x$ and velocity $\dot{x}$ of VP by solving (\ref{system})\\
  7: ~~~~~~~ $k=k+1$ \\
  8: \textbf{end while} \\
  \bottomrule
 \end{tabular}
 \end{center}
\end{table}

To prove stability of the optimal control algorithm in Table \ref{table_mpc}, we focus on proving the boundedness of the position error between the reference input and the output of the cognitive architecture.
In particular, optimality of the cost function is guaranteed in each optimization interval if the damped harmonic oscillator is adopted as end effector model instead of the HKB oscillator. Since both the reference position $r_p$ and the desired velocity $r_{\sigma}$ are bounded, we assume $r_p\in[\underline{r},\bar{r}]$ and $r_{\sigma}\in[\underline{v},\bar{v}]$.

\begin{theorem}
The optimal control algorithm applied to the HKB oscillator (\ref{system}) with cost function (\ref{vcost}) ensures bounded position error between the HP and the VP.
\end{theorem}

\begin{figure}
\scalebox{0.07}[0.07]{\includegraphics{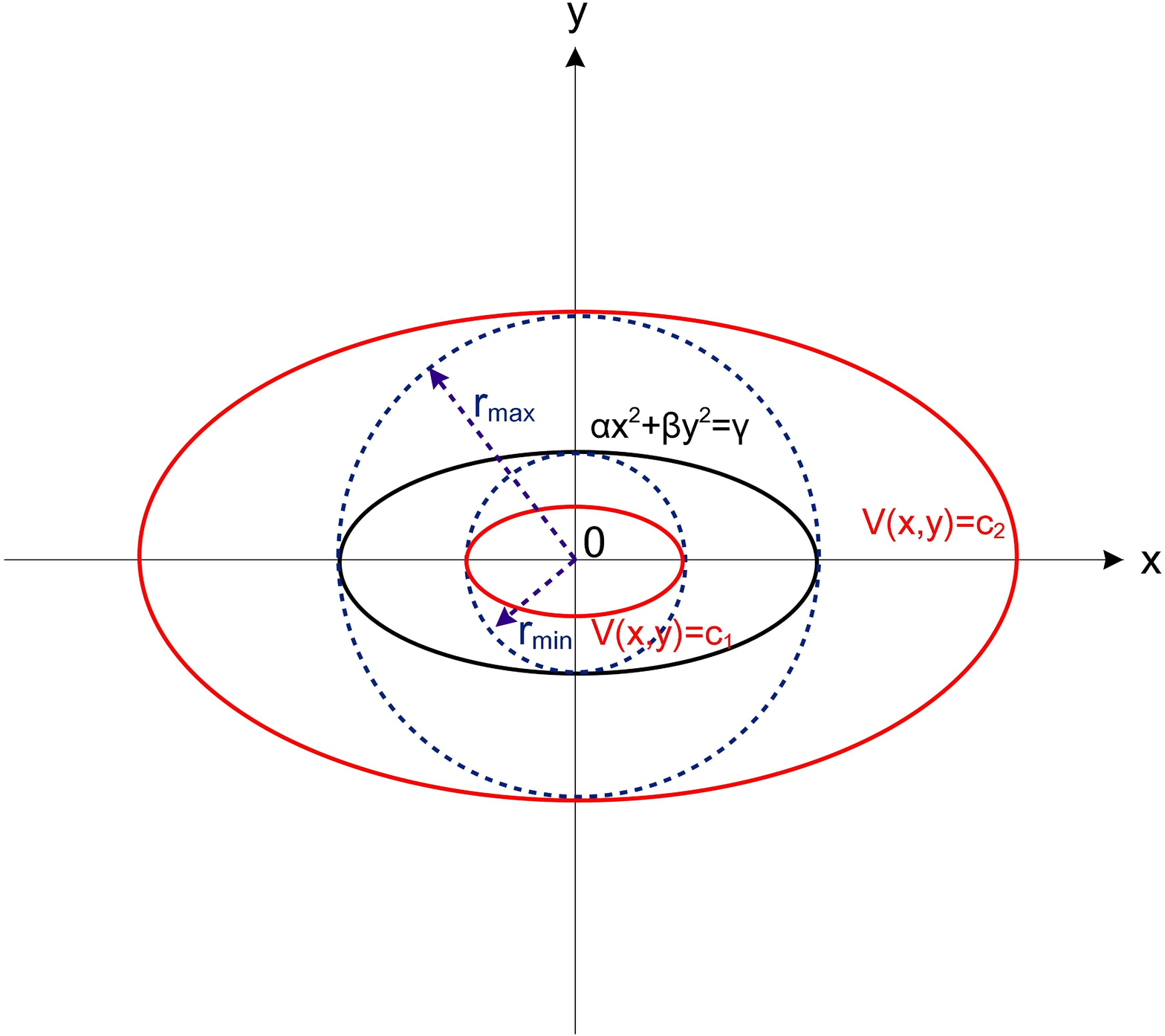}}\centering
\caption{\label{proof} Illustration on the construction of region $\mathfrak{R}$. The black ellipse is described by the equation $\alpha x^2+\beta y^2=\gamma$, and the region $\mathfrak{R}$ refers to the ring-shaped area bounded by two red ellipses corresponding to $V(x,y)=c_1$ and $V(x,y)=c_2$, respectively.}\centering
\end{figure}

Proof:
First of all, we need to demonstrate that there exists a limit cycle in the HKB oscillator
\begin{equation}\label{hkbo2}
\left\{
  \begin{array}{ll}
    \dot{x}=y \\
    \dot{y}=-(\alpha x^2+\beta y^2-\gamma)y-\omega^2 x
  \end{array}
\right.
\end{equation}
Choose the energy-like function as follows
$$
V(x,y)=\frac{\omega^2x^2+y^2}{2}
$$
The time derivative of $V(x,y)$ along the trajectory of the HKB oscillator (\ref{hkbo2}) is given by
\begin{equation*}
\begin{split}
\dot{V}(x,y)&=\omega^2x\dot{x}+y\dot{y}\\
&=\omega^2xy-(\alpha x^2+\beta y^2-\gamma)y^2-\omega^2 xy\\
&=-(\alpha x^2+\beta y^2-\gamma)y^2
\end{split}
\end{equation*}
Define
$$
r_{max}:=\max \left(\sqrt{\frac{\gamma}{\alpha}},\sqrt{\frac{\gamma}{\beta}} \right), r_{min}:=\min \left(\sqrt{\frac{\gamma}{\alpha}},\sqrt{\frac{\gamma}{\beta}} \right)
$$
and construct a region $\mathfrak{R}$ as follows (see Figure \ref{proof})
$$
\mathfrak{R}:=\{(x,y)\in \mathbb{R}^2: c_1 \leq V(x,y) \leq c_2\}
$$
where the positive constants $c_1$ and $c_2$ satisfy
\begin{equation*}\label{rmax}
r_{min} = \max\left(\sqrt{\frac{2c_1}{\omega^2}},\sqrt{2c_1} \right),\
r_{max} = \min\left(\sqrt{\frac{2c_2}{\omega^2}},\sqrt{2c_2} \right)
\end{equation*}
Clearly, $\mathfrak{R}$ contains no stationary points of the system. Indeed, the only stationary point of the system is $(x,y)=(0,0)$, but this stationary point is located outside of the region $\mathfrak{R}$. Moreover, $\dot{V}(x,y)\geq0$ when $V(x,y)=c_1$ and $\dot{V}(x,y)\leq0$ when $V(x,y)=c_2$. According to the Poincare-Bendixson theorem, we can conclude that the HKB oscillator (\ref{hkbo2}) has a limit cycle in $\mathfrak{R}$.

Let $J^*$ denote the value of the cost function (\ref{vcost}) with the optimal control algorithm in each time interval, and let $J_0$ represent the value of the corresponding cost function when $u=0$. Since $u$ aims at minimizing the value of the cost function for all $k\in \mathbb{N}^*$, we can write
\begin{equation*}
\begin{split}
J^*\leq J_0&=\frac{\theta_p}{2}(x(t_{k+1})-\hat{r}_p(t_{k+1}))^2 \\
&+\frac{\theta_{\sigma}}{2}\int_{t_k}^{t_{k+1}}(\dot{x}(\tau)-r_{\sigma}(\tau))^2d\tau
\end{split}
\end{equation*}
Recall that $r_p$ is bounded, which indicates that $\hat{r}_p(t_{k+1})$ is bounded according to equation (\ref{rvest}) and (\ref{rpest}). Moreover, $r_{\sigma}(\tau)$ is bounded as well, and note that $x(t)$ and $\dot{x}(t)$ are also bounded since the trajectory of the HKB oscillator converges to the limit cycle in $\mathfrak{R}$.
Thus, we can claim that $J_0$ is bounded for $k\in \mathbb{N}^*$, which implies boundedness of $J^{*}$ and as a consequence of the position error between the VP and the HP.
$\hfill{\blacksquare}$

\begin{remark}
It is demonstrated that the bound on the tracking error $|x(t_{k+1})-\hat{r}_p(t_{k+1})|$ converges to $0$ as $\theta_p\rightarrow 1$ and $\eta_m\rightarrow 0$. Similarly, the bound of the velocity error $|\dot{x}(t_{k+1})-r_{\sigma}(t_{k+1})|$ goes to $0$ if $\theta_{\sigma}\rightarrow1$, $\eta_m\rightarrow 0$ and $r_{\sigma}(t_k)=y(t_k)$ (see \emph{Supplementary Material} for the detailed analysis).
\end{remark}

The analytical solution for the optimization problem (\ref{minj}) is available if a linear damped harmonic oscillator is employed as the end effector model, and Pontryagin's minimum principle provides necessary and sufficient conditions to solve the minimization problem \cite{nai02}.

\begin{coro}\label{propv}
Given the linear system
$$
\ddot{x}+a\dot{x}+bx=u
$$
the optimal control approach guarantees convergence to the optimum solution over each subinterval.
\end{coro}

Proof: According to the fundamental theorem of the calculus of variations, we need to examine the second variation of the given cost function in order to establish the optimum. From the conclusions in \cite{nai02}, the second variation of the cost function (\ref{vcost}) is given by
\begin{equation*}
\begin{split}
\delta^2J&=\theta_p[\delta x(t_{k+1})]^2\\
&+\int_{t_k}^{t_{k+1}}\left(
                                  \begin{array}{cc}
                                    \delta X & \delta u \\
                                  \end{array}
                                \right)\left(
                                         \begin{array}{cc}
                                           H_{XX} & H_{Xu} \\
                                           H^T_{uX} & H_{uu} \\
                                         \end{array}
                                       \right)\left(
                                                \begin{array}{c}
                                                  \delta X \\
                                                  \delta u \\
                                                \end{array}
                                              \right)dt
\end{split}
\end{equation*}
where $H$ is the Hamiltonian
$$
H(X,u,\lambda)=\frac{1}{2}\theta_s(\dot{x}-r_{\sigma})^2+\frac{1}{2}\eta u^2+\lambda^T\left(
          \begin{array}{c}
            y \\
            -ay-bx+u \\
          \end{array}
        \right)
$$
with $X=[x,\dot{x}]^T=[x,y]^T$ and $\lambda=[\lambda_1,\lambda_2]^T$.
Rewrite the linear system in matrix form as follows
$$
\dot{X}=AX+Bu
$$
where
$$
A=\left(
    \begin{array}{cc}
      0 & 1 \\
      -b & -a \\
    \end{array}
  \right)
, \quad B=\left(
                         \begin{array}{c}
                           0 \\
                           1 \\
                         \end{array}
                       \right)
$$
Let $X=X^*+\delta X$ and $u=u^*+\delta u$, where $X^*$ and $u^*$ denote optimal state and optimal control, respectively. Since $\dot{X}^*=AX^*+Bu^*$, we get
\begin{equation}\label{delX}
\dot{\delta X}=A \delta X+B\delta u
\end{equation}
where $\delta X=[\delta x, \delta\dot{x}]^T$. Thus, it follows from $H_{Xu}=H_{uX}=[0~0]^T$, $H_{uu}=\eta>0$ and
$$
H_{XX}=\left(
         \begin{array}{cc}
           0 & 0 \\
           0 & \theta_s \\
         \end{array}
       \right) \geq 0
$$
that
\begin{equation*}
\begin{split}
\delta^2J&=\theta_p[\delta x(t_{k+1})]^2+\int_{t_k}^{t_{k+1}}\delta X(t)^TH_{XX}\delta X(t)+\eta(\delta u(t))^2dt\\
&=\theta_p[\delta x(t_{k+1})]^2+\int_{t_k}^{t_{k+1}}\theta_s(\delta \dot{x}(t))^2+\eta(\delta u(t))^2dt\\
&\geq0
\end{split}
\end{equation*}
Moreover, $\delta^2J=0$ is equivalent to $\delta x(t_{k+1})=0$, $\delta \dot{x}(t)=0$ and $\delta u(t)=0$ for all $t\in [t_k,t_{k+1}]$, which yields $\delta x(t)=\delta x(t_k)=0$ from (\ref{delX}). This corresponds to the optimal solution $X^*$ and the optimal control $u^*$. Therefore, we conclude that the optimal control ensures the minimum value of the cost function (\ref{vcost}) for the linear system in each time interval. \hfill $\blacksquare$

\section{Validation}\label{sec:validation}

In this section we experimentally validate our control algorithms on a simple, yet effective, set-up implemented at the University of Bristol, UK. Experimental data of human-human interaction is used to evaluate the matching performance of the VP and a comparison with existing VP models is provided as well.

\subsection{Experimental Set-up}

\begin{figure}
 \centering
 \includegraphics[width=0.6\textwidth]{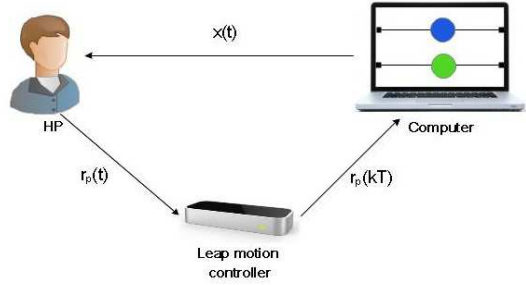}
\caption{Experimental set-up of the mirror game between a HP and a VP. The position of the human fingertip $r_p(t)$ is detected by a leap motion controller, and the sampled position $r_p(kT)$ is sent to the computer, while the position $x(t)$ of the VP is generated by implementing the control algorithm. Two balls are shown on the computer screen, which describe the end effector positions of the HP (green ball) and the VP (blue ball), respectively.}\label{fig:ExSetup}
\end{figure}

The employed set-up was developed for measuring motions of players in the one-dimensional mirror game. A human participant is required to join the game while interacting with a VP (implemented on
a laptop computer). In order to detect the position of his/her hand, a leap motion controller \cite{leap} is employed (see Fig. \ref{fig:ExSetup}). The leap motion controller and the laptop computer are both placed on a table whose height is around 70cm. The HP is required to wave his/her hand horizontally over the leap motion controller at a vertical distance of approximately 50cm. Indeed, at this distance the horizontal resolution of the device is maximum and it is able to measure the hand position within a range of 60cm. The position of the hand of the HP within this interval is mapped intto the interval $[-0.5,0.5]$ and visualized on the computer screen as a green solid circle, while the position of the VP is visualized as a blue solid circle. The adaptive control algorithm described in Section \ref{sec:adaptive} is implemented with Euler method in Matlab (version R2012b), whilst the solver ``bvp4c'' is adopted to handle the optimization problem of the optimal control algorithm presented in Section \ref{sec:optimal}. Players can be either standing or seated. After the game is initialized, there is a 2s wait before data recording begins and the game starts. This initial delay is used to allow the HP to place his/her hand over the leap motion controller. Human players are not instructed before playing the game, but they are just told to act as a leader and let the VP follow them during a 60s round. The case where the VP acts as a leader is also tested experimentally.

\subsection{Measures}
The temporal correspondence between the VP and the HP is evaluated according to the following indexes:
the root mean square (RMS), the relative position error (RPE)\cite{piotr15}, the circular variance (CV)\cite{kreuz07} and the time lag (TL)\cite{orfan96}.
\begin{enumerate}
  \item RMS: The root mean square of the position error between two players describes the tracking accuracy of the follower in the mirror game.
  $$
  RMS=\sqrt{\frac{1}{n}\sum_{i=1}^{n}(x_{1,i}-x_{2,i})^2}
  $$
  where $n$ is the number of sampling steps in the simulation, and $x_{1,i}$ and $x_{2,i}$ denote the positions of the leader and the follower at the $i$-th sampling step, respectively.
  \item RPE: The relative position error is a measure of how well the follower was tracking the leader in the mirror game \cite{piotr15}. Positive values of the RPE indicate that the follower is indeed behind the leader.
  $$
  RPE=\left\{
        \begin{array}{ll}
          &(x_1(t)-x_2(t))sgn(\dot{x}_1(t)),\\
          &~~~~~~~~~~~~~\textmd{if}~\hbox{$sgn(\dot{x}_1(t))=sgn(\dot{x}_2(t))\neq0$;}\\
          &|x_1(t)-x_2(t)|,\hbox{otherwise.}
          %\\&~~~~~~~~~~~~~~\hbox{otherwise.}
        \end{array}
      \right.
  $$
  where $x_1(t)$ and $x_2(t)$ ($\dot{x}_1(t)$ and $\dot{x}_2(t)$) are the positions (velocities) of the leader and follower at time $t$, respectively.
  \item CV: The circular variance is used to quantify the coordination level between two players
  $$
  CV=\left\|\frac{1}{n}\sum_{k=1}^{n}e^{i\Delta \Phi_k}\right\| \quad \in [0,1]
  $$
  where $\Delta \Phi_k$ represents the relative phase between two players at the $k$-th sampling step, $n$ refers to the total number of time steps and $\|\cdot\|$ denotes the $2$-norm.
  \item TL: The time lag describes the shifted time that achieves the maximum cross-covariance of the two time series. It is sensitive to the changes of motion direction and hence can be interpreted as the average reaction time of the player in the mirror game \cite{orfan96}.
\end{enumerate}

\subsection{Results}

\subsubsection{VP driven by the ICA based on adaptive control}

The parameters for the HKB equation and the adaptive feedback controller (AFC) in (\ref{eq:u}) are set heuristically as follows:
$\alpha=10$, $\beta=20$, $\gamma=-1$, $\omega=0.1$, $a(0)=-5$, $b(0)=-5$, $C_p=40$ and $\delta=0.25$. In our implementation the sampling time is $T=0.1$s and therefore $\eta_a=30>\frac{\ln 2}{2T} \simeq 0.35$. In particular, the values of all the previous parameters have been chosen so that the response of the HKB oscillator to several sinusoidal signals with different frequencies is qualitatively the same as the one of a HP trying to track the same references. Moreover, it is worth pointing out that the initial values of $a$ and $b$ influence the performance of the avatar only at start-up.

The reactive-predictive controller (RPC) proposed in \cite{noy11} is also implemented to compare its performance against that of our adaptive feedback controller when considering the same input trajectory from the human leader.  Following the scheme presented in \cite{noy11} to implement the RPC, the dynamics of the VP is described by the following system:
\begin{equation*}
\ddot{x} = \sum_{i=1}^5 A_i \omega_i cos (\omega_i t) + f
\end{equation*}
where $x \in \mathbb{R}$ represents the position of the avatar and
\begin{equation*}
\dot{f} = k(\hat{r}_v-\dot{x}), \quad k >0
\end{equation*}
with the parameters $A_i$ being estimated adaptively as
\begin{equation*}
\dot{A}_i = \lambda \left [ \hat{r}_v-\sum_{i=1}^5 A_i sin (\omega_i t) \right ] sin (\omega_i t), \quad \lambda >0
\end{equation*}
As suggested in \cite{noy11}, in this case the parameters are chosen as follows:
$\omega_1 = 0.025$, $\omega_2 = 0.05$, $\omega_3 = 0.075$, $\omega_4 = 0.1$, $\omega_5 = 0.125$, $\lambda = 0.01$, $k = 30$ and $A_i(0) = 0, \forall i =1,...,5$.

\begin{figure}\centering
 \subfigure[position time series]
 {\includegraphics[width=0.45\textwidth]{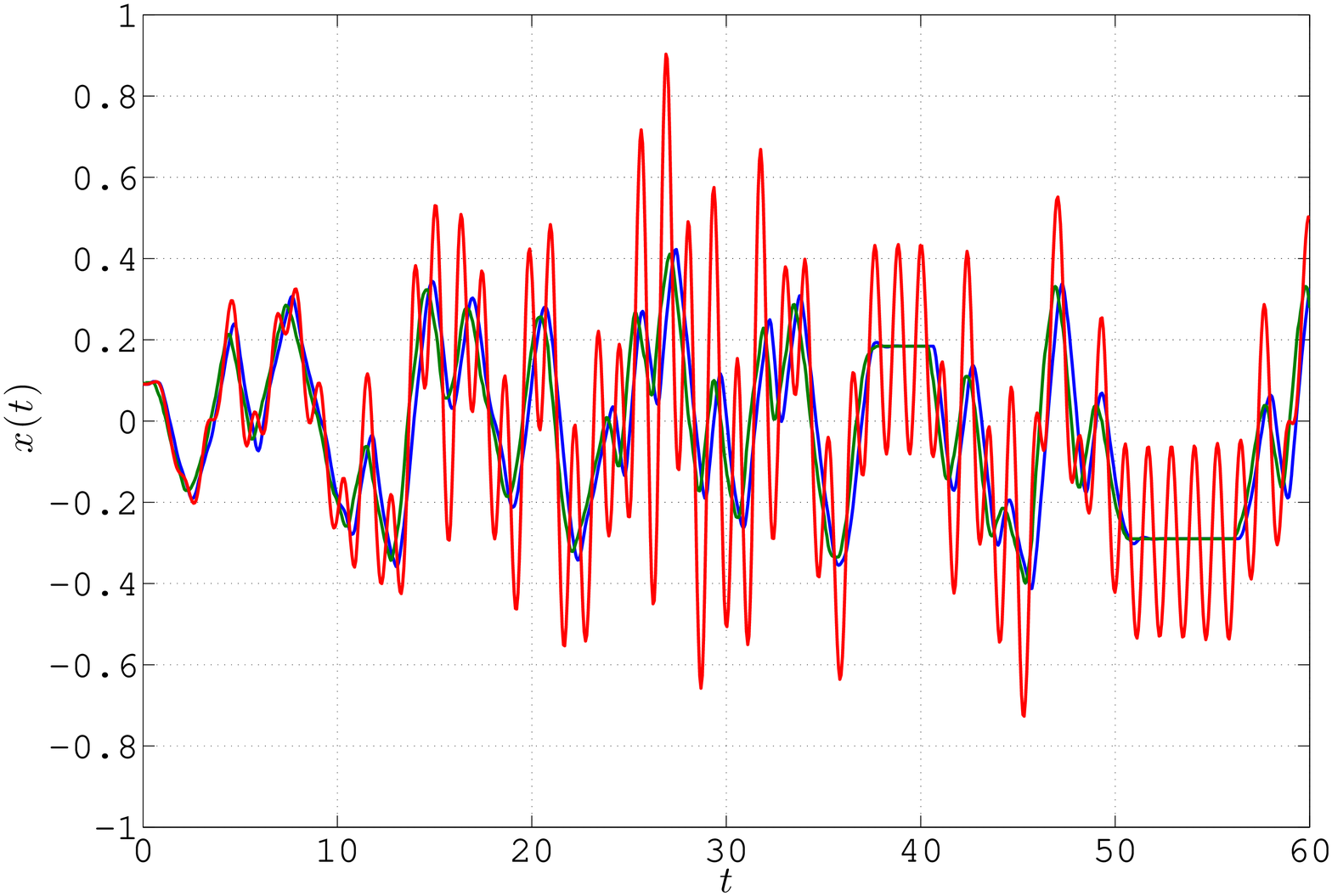}}
 \subfigure[relative phase time series]
 {\includegraphics[width=0.45\textwidth]{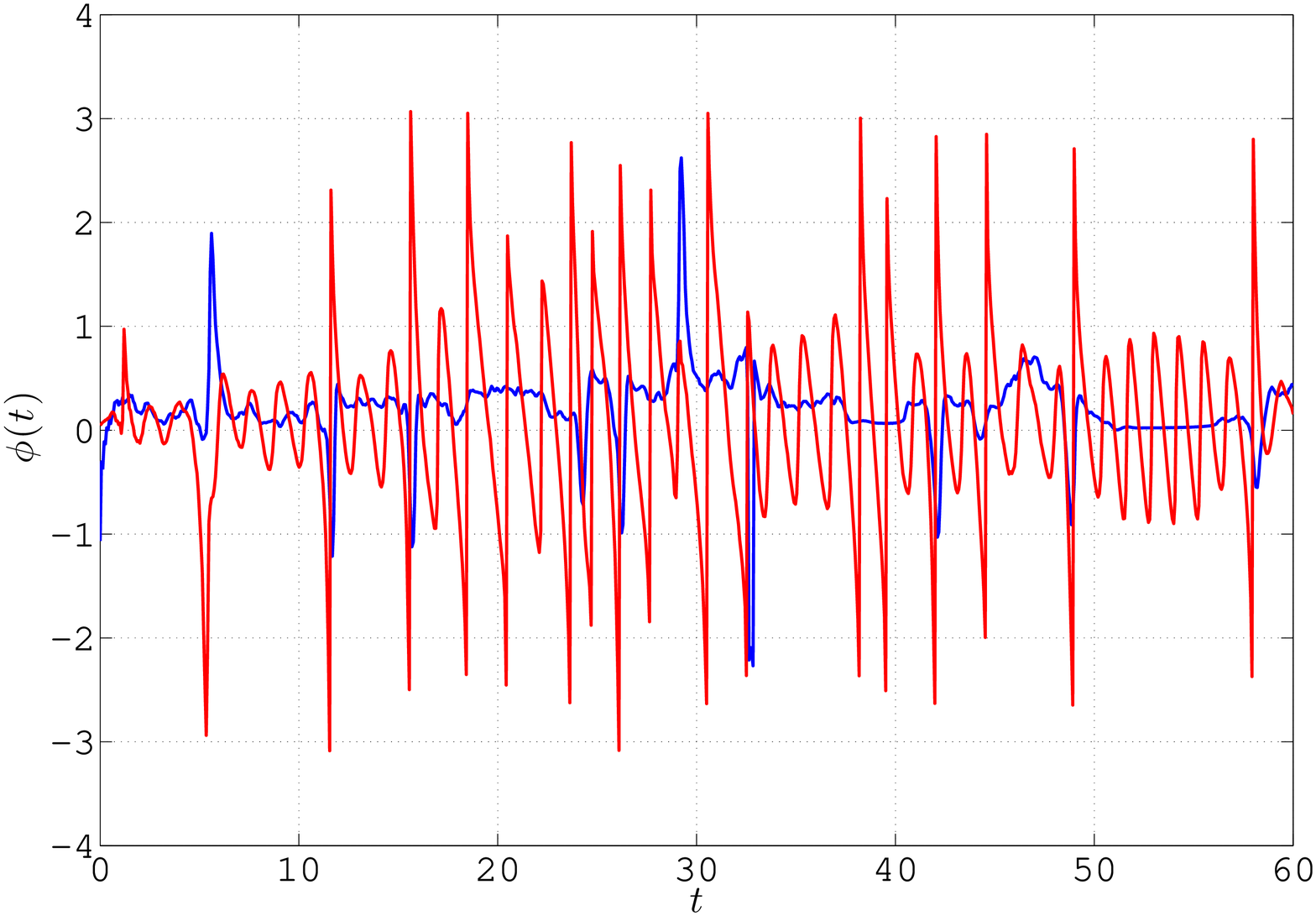}}
 \caption{Time series of the position (a) and of the relative phase (b) between the human leader and the VP; blue (AFC), red (RPC), green (human leader)}
 \label{fig:tspics}
\end{figure}

\begin{figure}\centering
 \subfigure[veocity PDF]
 {\includegraphics[width=0.45\textwidth]{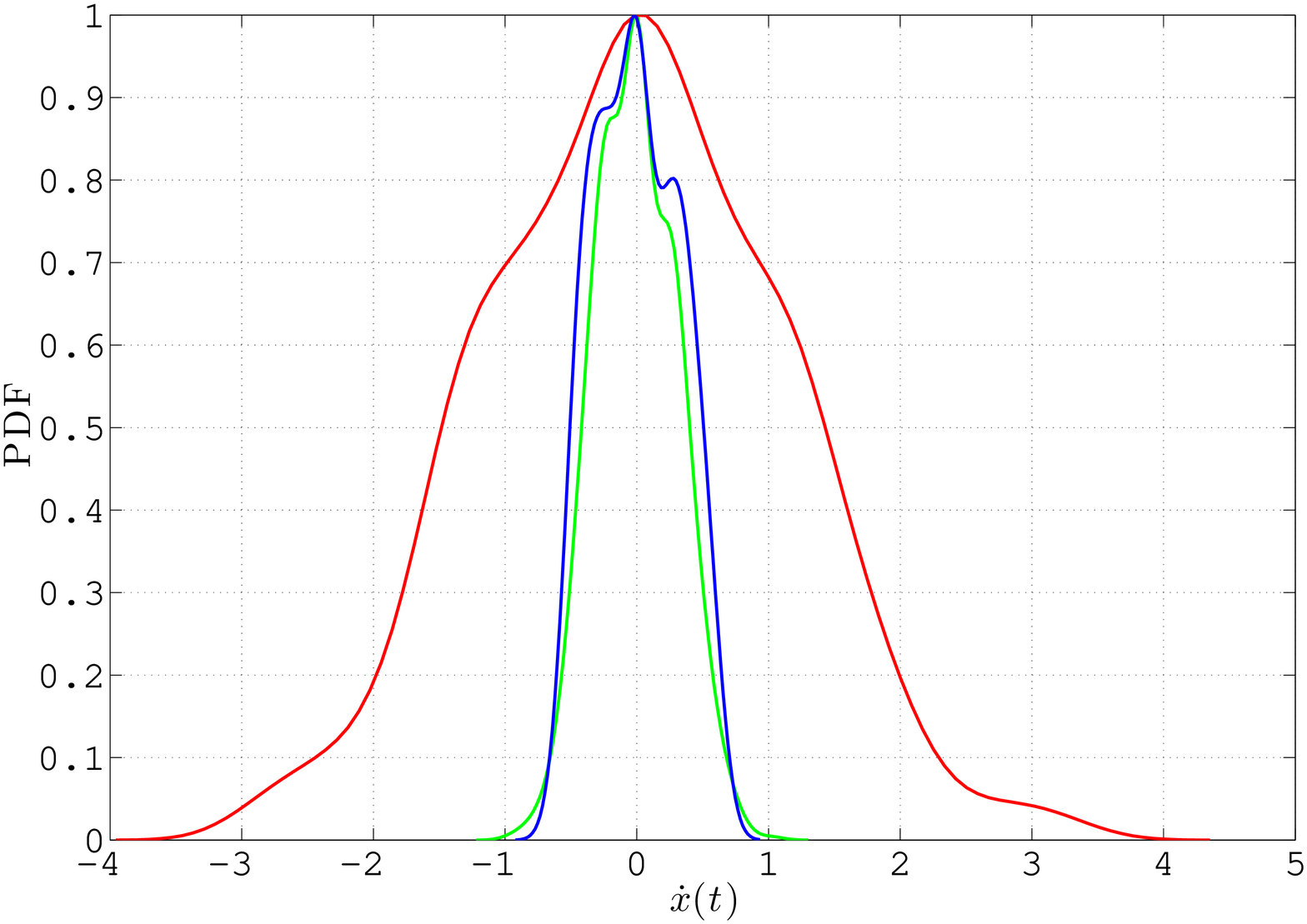}}
 \subfigure[relative phase PDF]
 {\includegraphics[width=0.45\textwidth]{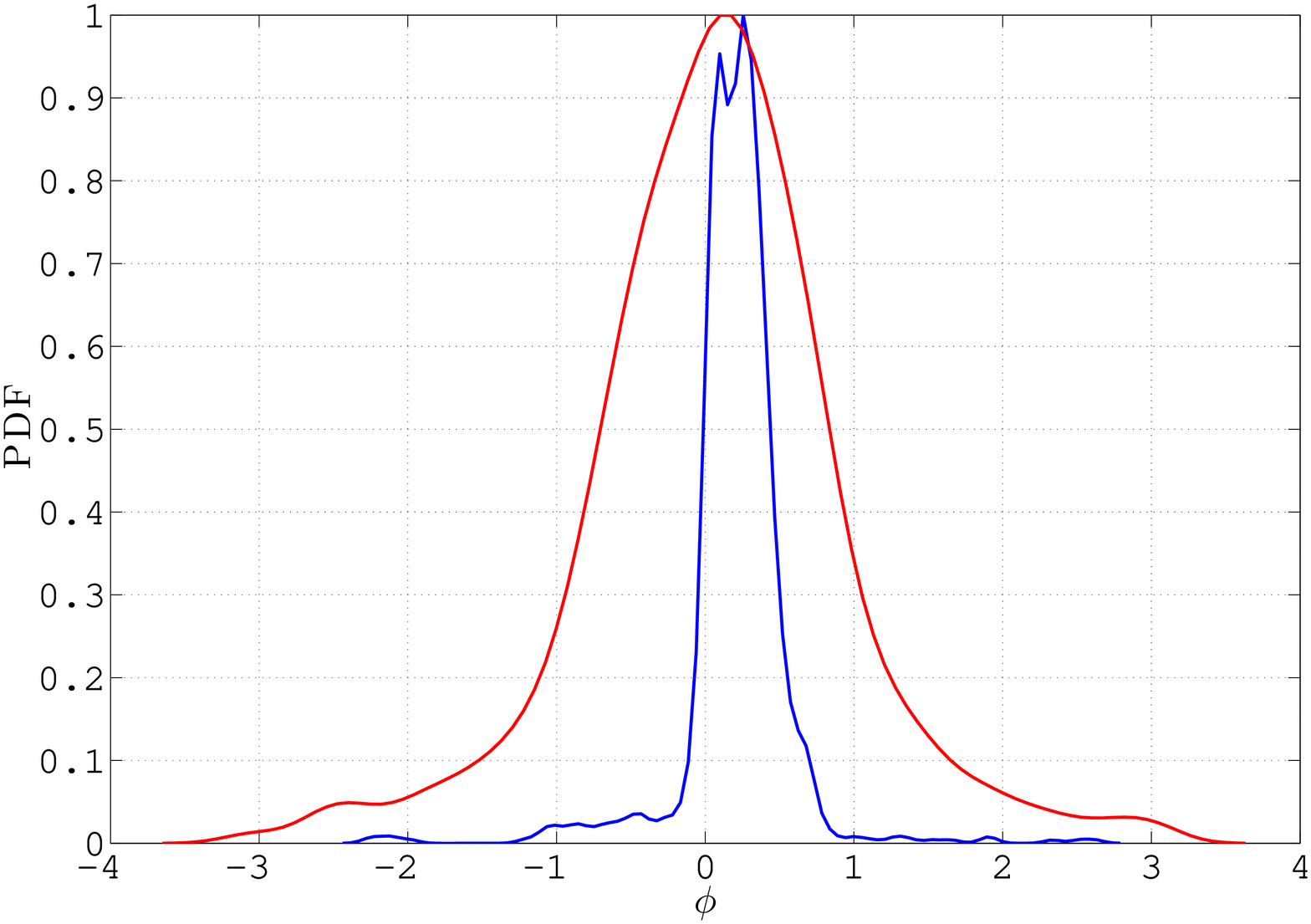}}
 \caption{Distributions of the velocity (a) and the relative phase (b) between the human leader and the VP; blue (AFC), red (RPC), green (human leader)}
 \label{fig:pdfpics}
\end{figure}

To compare the performance of the two algorithms, we plot the time series of both the position and the relative phase (see Fig. \ref{fig:tspics}) together with the distributions of velocity and relative phase of the human and the VP (see Fig. \ref{fig:pdfpics}). In particular, the relative phase between the two players is defined as $\Delta \phi=\phi_{HP}-\phi_{VP}$, where $\phi_{HP}$ and $\phi_{VP}$ are the phases of the HP and the VP, respectively. In addition, the phase is estimated according to the method proposed in \cite{kral08}. Note that positive values of $\Delta \phi$ correspond to the avatar following the HP during the game.

We can observe that, when using the reactive-predictive controller, the position of the VP presents oscillations away from the human participant position not only when he/she is moving, but also when he/she is still. Such an oscillatory feature does not appear when using the adaptive feedback controller. In general, both the position error $e=x-r_p$ and the velocity error $\dot{e}=\dot{x}-\hat{r}_v$ turn out to be higher when using the RPC strategy.  When using the adaptive feedback controller, the position error remains smaller never exceeding $0.2$, while it can become as high as $0.8$ when using RPC. Similarly, the velocity error never exceeds $0.62$ for the AFC, while it goes up to a maximum of $3$ for the RPC.

Moreover, when using AFC, the relative phase time series is much closer to $0$ than that obtained when using RPC, meaning that with our proposed algorithm it is possible for the VP to better synchronize with the human leader. Such results are confirmed by the relative phase distributions obtained when using both the algorithms, as shown in Fig. \ref{fig:pdfpics}(b).
Finally, the difference in the velocity distributions of the HP and the VP is much more evident when RPC is used, confirming that our strategy better captures the features of the HP and is therefore able to replicate more accurately the kinematic properties observed in human motor coordination in the context of the mirror game.

\subsubsection{VP driven by the ICA based on optimal control} \label{MPC}

\begin{figure}\centering
 \subfigure[Temporal Correspondence]{\includegraphics[width=0.4\textwidth]{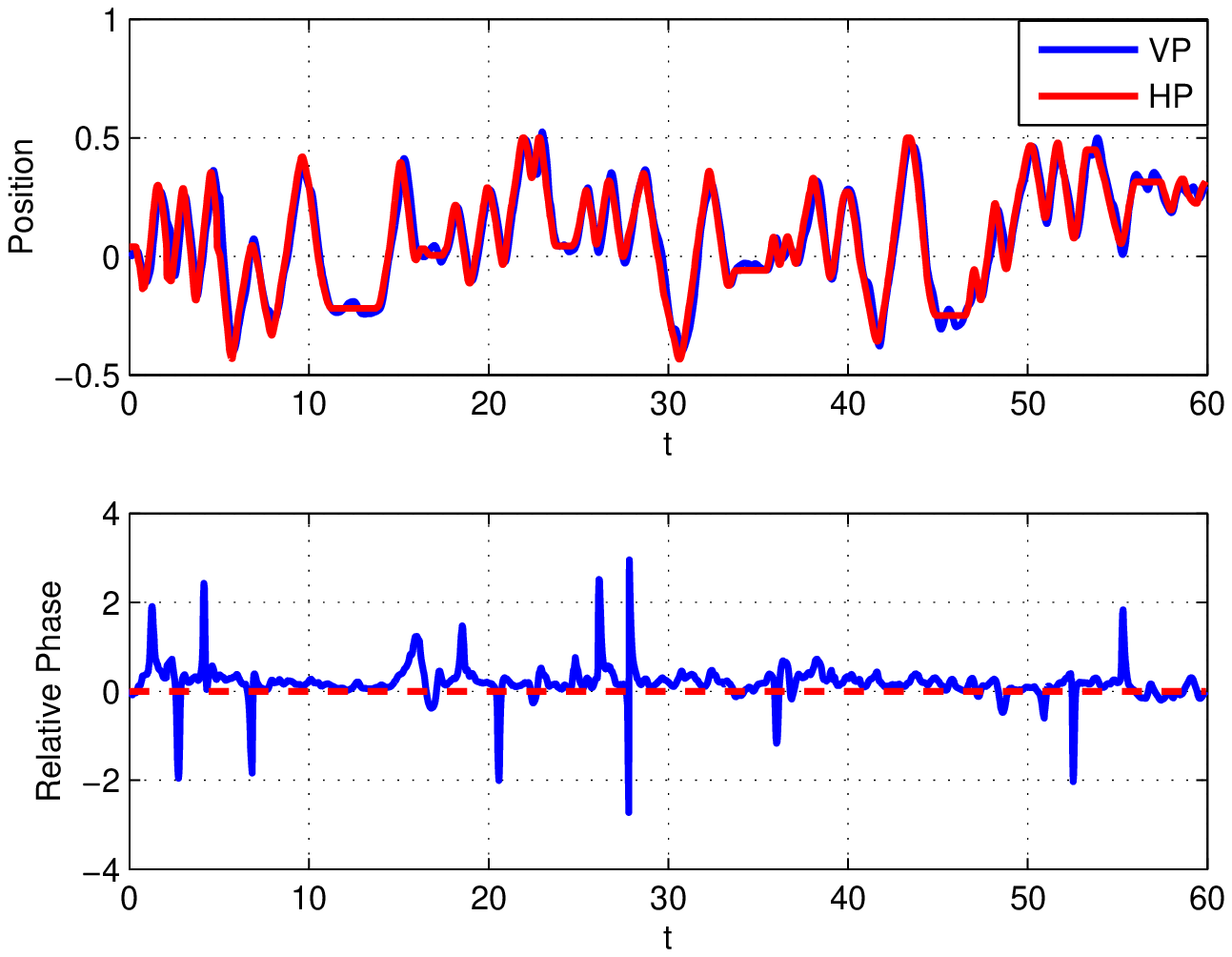}}
 \subfigure[Similarity]{\includegraphics[width=0.4\textwidth]{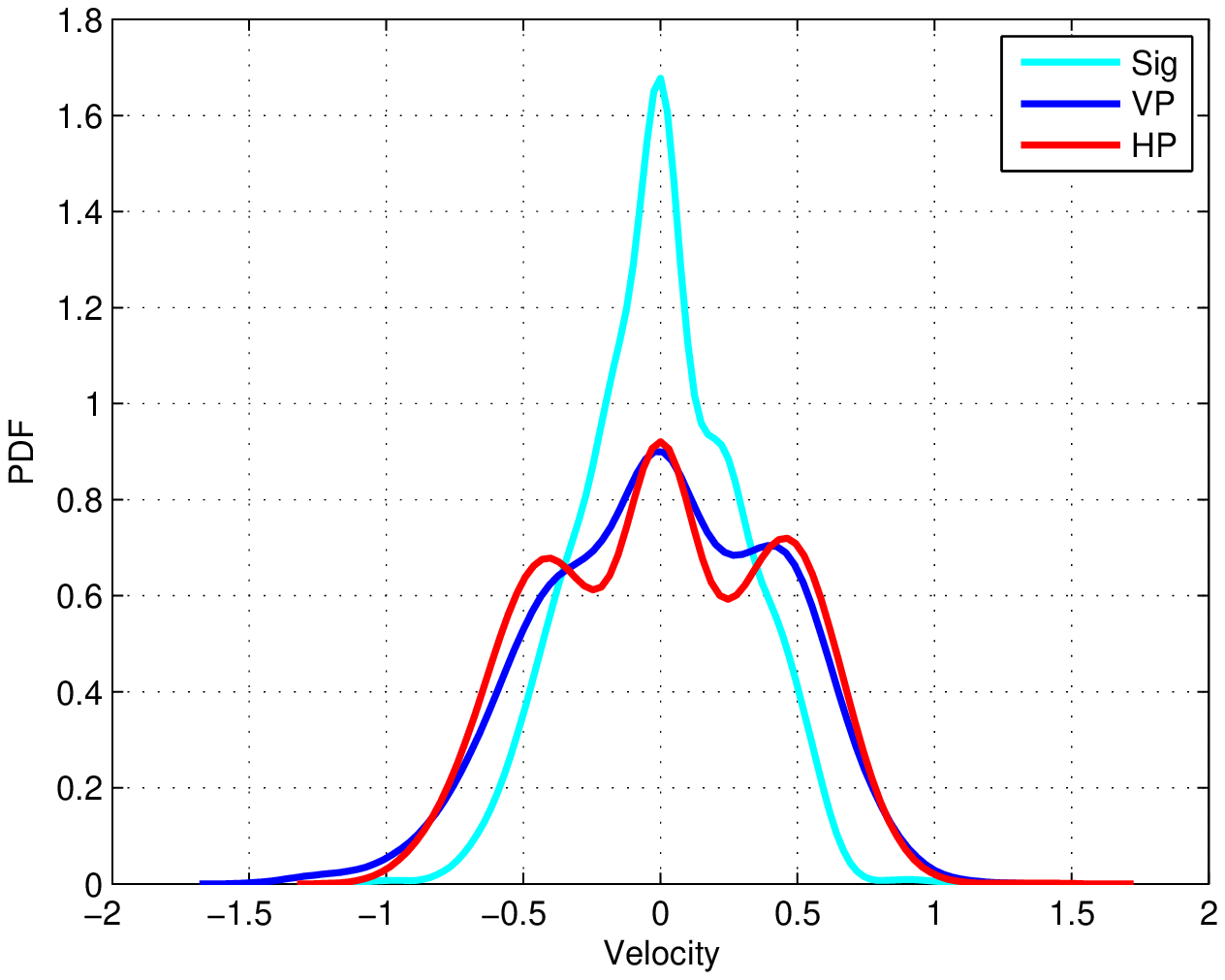}}
 \caption{Time evolution of positions and relative phase (a) and PDF of velocities (b) while the VP is driven by the optimal control and acts as follower in the mirror game}\label{fig:follower}
\end{figure}

\begin{figure}\centering
 \subfigure[Temporal Correspondence]{\includegraphics[width=0.4\textwidth]{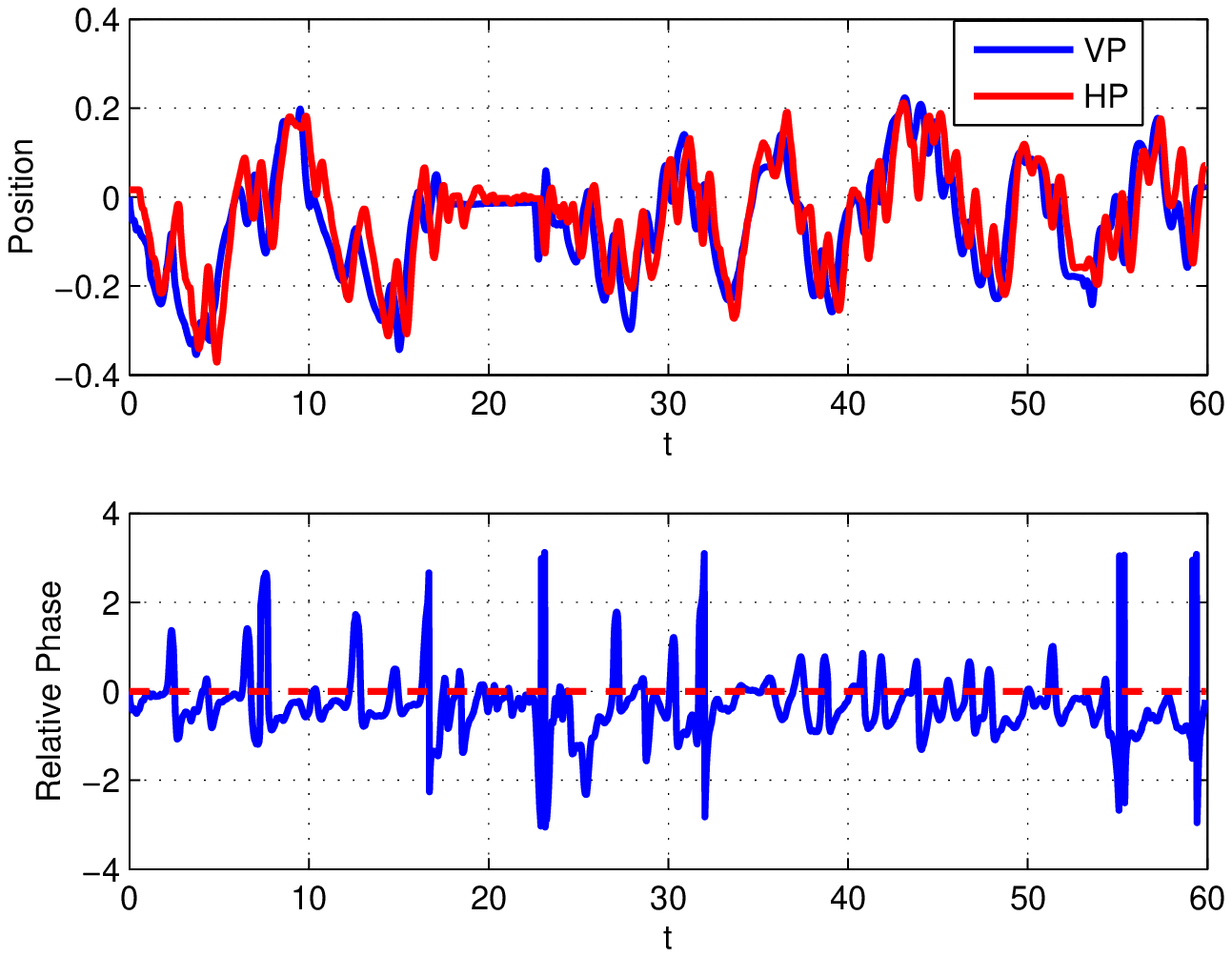}}
 \subfigure[Similarity]{\includegraphics[width=0.4\textwidth]{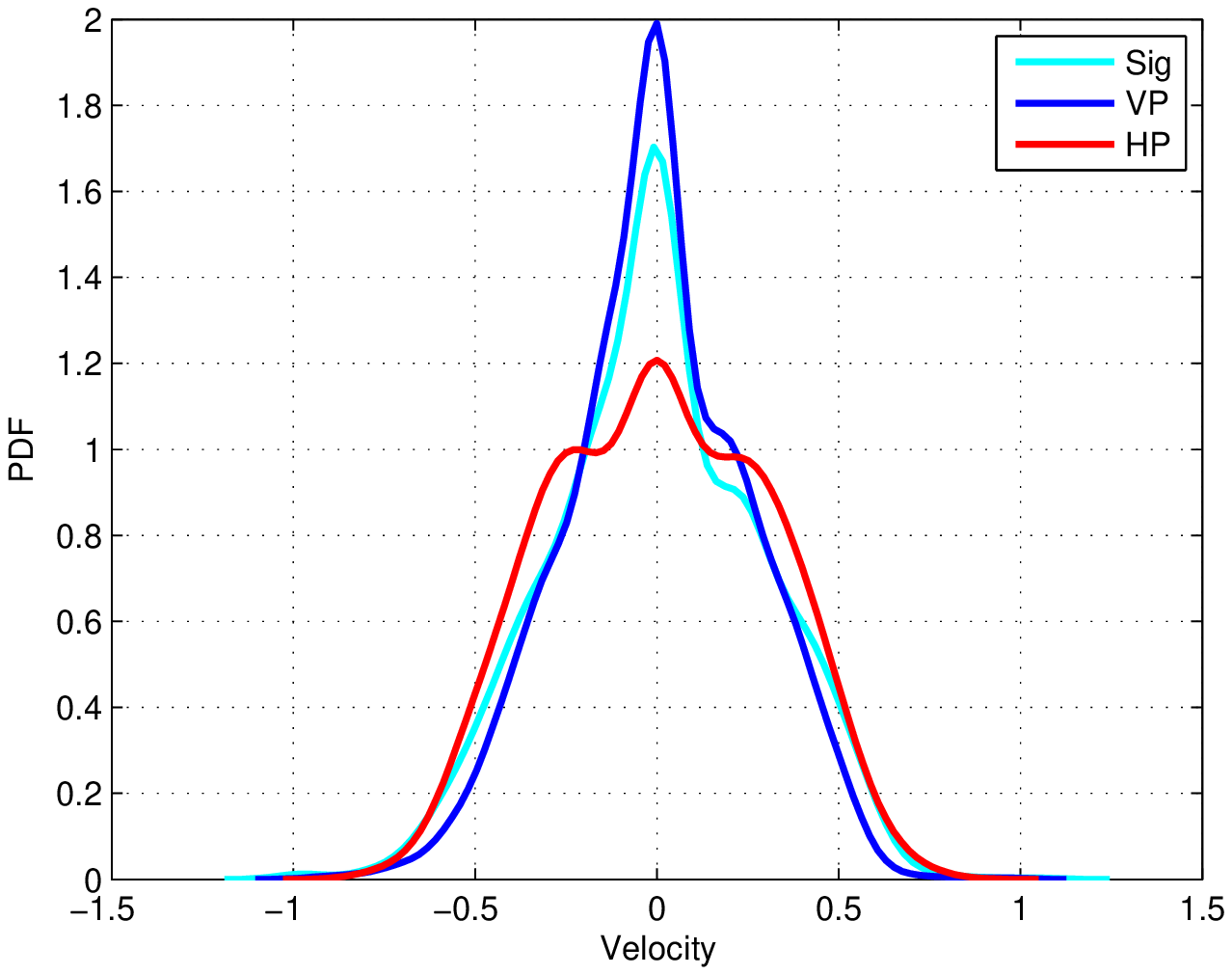}}
 \caption{Time evolution of positions and relative phase (a) and PDF of velocities (b) while the VP is driven by the optimal control and acts as leader in the mirror game}\label{fig:leader}
\end{figure}

The parameters of the VP are set heuristically as follows: $\alpha=1$, $\beta=1$, $\gamma=1$, $\omega=1$, $\eta_m=10^{-4}$ and the sampling period $T=t_{k+1}-t_k=0.03$s. In order for the VP to play the mirror game as a follower, we set the control parameters $\theta_p=0.9$ and $\theta_{\sigma}=0.1$, which makes the VP play in a follower configuration (minimizing the position mismatch more than the signature EMDs). As we can see from the top panel in Fig.\ref{fig:follower}(a), the VP performs quite well as a follower during the game; indeed, the root mean square (RMS) of the position error is equal to $0.057$. In order to distinguish the leader from the follower in the game, we also calculate the relative phase between the HP and the VP. From the bottom panel in Fig. \ref{fig:follower}(a) we can observe that the majority of relative phase is positive, meaning that the VP is following the HP in the game for most of the time. The circular variance (CV) is also calculated to take into account the coordination level between two players. The CV between the HP and the VP is $0.95$, which indicates a high coordination level. As for the distribution of the velocity, we can see in Fig. \ref{fig:follower}(b) that the VP signature (blue line) is closer to that of the HP (red line) than the desired motor signature (cyan line). This is due to the choice of the control parameters in the cost function (\ref{vcost}) that render the strategy able to minimize more the position error between the players. The measured EMDs at the end of the trial are given as follows: EMD$(Sig,HP)$ = EMD$(Sig,VP)=0.017$ and EMD$(VP,HP)=0.005$. \\

The VP can be enabled to play the game as a leader by changing the control parameters setting $\theta_p=0.1$ and $\theta_{\sigma}=0.9$. Experimental results are shown in Fig. \ref{fig:leader}. The RMS of the tracking error is $0.08$, and the CV between the two players is $0.81$. As depicted in the bottom panel of Fig. \ref{fig:leader}(a), the majority of the relative phase time series is negative, meaning that now the VP is leading the HP during the game for most of the time. In contrast to the previous case, the velocity distributions shown in Fig. \ref{fig:leader}(b) confirm that the VP is now matching well the desired signature (velocity profile). In this case the trade off is slightly larger but still the relevant EMDs are given as follows: EMD$(Sig,HP)$ =EMD$(Sig,VP)=0.004$ and EMD$(VP,HP)=0.008$.

\subsubsection{Interaction between two VPs}

\begin{figure}
 \centering
  \includegraphics[width=0.7\textwidth]{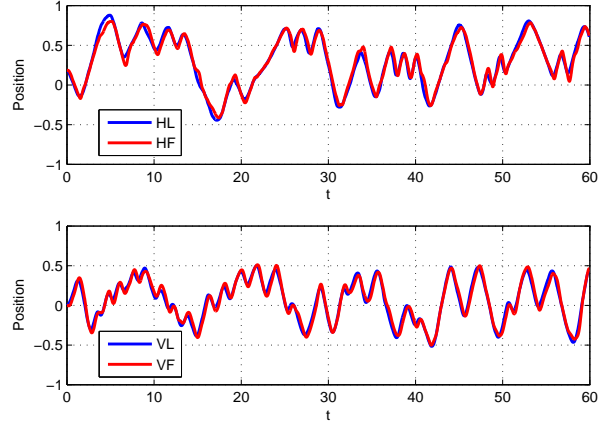}
  \caption{\label{pos_2vp}Position time series of the HP pair (upper panel) and the corresponding VP pair (lower panel) in the mirror game.}
\end{figure}

As mentioned before, our final goal is to create a customized VP able to ``replicate'' the kinematic features of a given HP in the mirror game. To test how well the VP can replicate the features of a given HP, we carried out the following experiment. First of all, two HPs are required to play the mirror game in a Leader-Follower condition. Then the signatures of the human leader (HL) and the human follower (HF) are fed into a virtual leader (VL) and a virtual follower (VF), respectively. Finally, we make the VL ($\theta_p=0.43$) and the VF ($\theta_p=0.92$) play the mirror game together. In Fig. \ref{pos_2vp}, the upper panel shows the time evolution of the position trajectories of the HL and the HF, while the lower panel presents the position trajectories of the corresponding VL and VF. It appears that the VL and the VF succeed in matching the kinematic characteristics of the HL and the HF in terms of the RMS value of their position error [RMS$(HL,HF)=0.0466$ and RMS$(VL,VF)=0.0497$] and the time lag between the two players [TL$(HL,HF)=0.09$ and TL$(VL,VF)=0.09$].
In addition, Table \ref{table:emd} gives the values of EMDs and describes the matching results quantitatively. $Sig_L$ and $Sig_F$ represent the signatures of HL and HF when playing solo, respectively. It shows that the proposed approach allows to replicate effectively the game dynamics between two humans playing the mirror game via two coupled VPs.

\begin{table}
\caption{\label{table:emd} Matching performance of VPs in terms of EMD.}
\begin{center}
\begin{tabular}{|c|c|c|c|}
  \hline
  EMD($Sig_{L}$,HL) & 0.010 & EMD($Sig_{F}$,HF) & 0.007    \\ \hline
  EMD($Sig_{L}$,VL) & 0.006 & EMD($Sig_{F}$,VF) & 0.006    \\ \hline
  EMD(HL,HF)  & 0.0034 & EMD(VL,HL) & 0.0052   \\ \hline
  EMD(VL,VF)  & 0.0031 & EMD(VF,HF) & 0.0053   \\ \hline
\end{tabular}
\end{center}
\end{table}

\subsubsection{Comparison with existing models}

\begin{figure}
\scalebox{2}[2]{\includegraphics[width=1.5in]{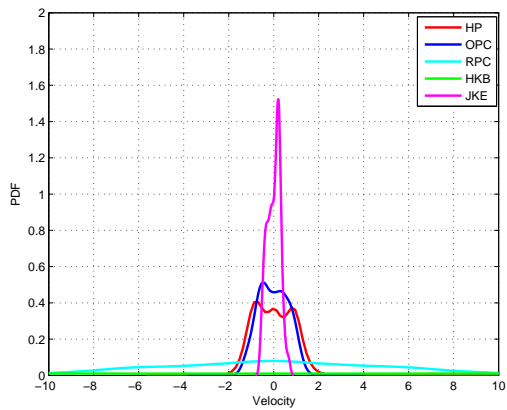}}\centering
\caption{\label{vel_dis}PDF of velocity time series for different VP models.}
\end{figure}

In order to compare the VP models, we need to establish a benchmark, which describes the general kinematic characteristics of human participants in the mirror game. To this aim, $5$ human participants were asked to track a prerecorded reference signal, and indexes of temporal correspondence were recorded to represent a benchmark of typical human dynamics. The existing VP models were then enabled to track the same reference. Corresponding indexes were computed for the VPs and compared with the benchmark to evaluate the proposed control algorithms. We adopted the following models to drive the VP: our strategy based on optimal control (OPC), Haken-Kelso-Bunz model (HKB)\cite{kel09}, reactive-predictive control (RPC)\cite{noy11} and Jirsa-Kelso excitator (JKE) \cite{pnas14}. The parameter setting of existing VP models is the same as that in the relevant references. In Table \ref{table:comparison}, we show the benchmark of temporal correspondence and performance indexes of the corresponding VP models. On the whole, our algorithm performs best in terms of matching the benchmark among all the VP models. In addition to temporal correspondence, we also computed how similar the VP signature is to that of the HP during the benchmark experiment (see Fig.\ref{vel_dis}). We find again that the optimal control strategy developed in this paper is the best in terms of replicating the human-like movement with an estimated EMD$(HP,OPC)=0.0184$.

\begin{table}
\caption{\label{table:comparison} Indexes of temporal correspondence}
\begin{center}
\begin{tabular}{|c|c|c|c|c|c|}
  \hline
      & HP & OPC  & HKB & RPC & JKE \\ \hline
  RPE & 0.0914 & 0.0816  & 2.1767 & 0.3838 & 0.1467    \\ \hline
  CV  & 0.3011 & 0.1002  & 0.9400 & 0.7602 &  0.4859   \\ \hline
  TL  & 0.2035 & 0.1274  & 1.5192 & 0.0428 &  -0.9674  \\ \hline
\end{tabular}
\end{center}
\end{table}

\section{Conclusions}\label{sec:conclusions}
We presented the novel design of an interactive cognitive architecture able to drive a virtual player to play the mirror game against a human player. Two strategies were developed. The first, based on adaptive control, was shown to be effective to achieve temporal correspondence between the motion of the virtual player and that of the human individual. Convergence of the algorithm was proved. It was noticed that the adaptive control strategy does not allow the VP to exhibit some desired kinematic features (individual motor signature) of a given human player. To overcome this limitation, a different strategy based on the iterative solution of an appropriate optimal control problem was proposed. After proving boundedness and convergence of this additional approach, its effectiveness was tested experimentally. It was shown that the proposed strategy is able to drive the VP so as to play the game both as leader or follower while matching well the individual motor signature of a given individual. Finally, a comparison with other existing models was carried out confirming the effectiveness of the proposed approach. We wish to emphasize that our approach opens the possibility of making VPs, each modeling a different individual, play against each other and produce in silico experiments. This can reduce the cost and time of carrying out joint action experiments and can be effectively used to test different human-machine interaction scenarios via the mirror game.

\section*{Acknowledgements}
This work was funded by the European Project AlterEgo FP7 ICT 2.9 - Cognitive Sciences and Robotics, Grant Number 600610. The authors wish to thank Prof. Benoit Bardy, Dr. Ludovic Marin and Dr. Robin Salesse at EUROMOV, University of Montpellier, France for all the insightful discussions and for collecting some of the experimental data that is used to validate the approach presented in this paper.

\section*{Supplementary Material}
The solution of the optimal control algorithm in each time interval $[t_k,t_{k+1}]$ can be transformed into a boundary value problem by applying Pontryagin's minimum principle \cite{nai02}. We start by constructing the Hamiltonian as follows
$$
H(X,u,\lambda)=\frac{1}{2}\theta_{\sigma}(\dot{x}-r_{\sigma})^2+\frac{1}{2}\eta_m u^2+\lambda^T\left(
          \begin{array}{c}
            y \\
            -(\alpha x^2+\beta y^2-\gamma)y-\omega^2x+u \\
          \end{array}
        \right)
$$
where $X=[x,\dot{x}]^T=[x,y]^T$ and $\lambda=[\lambda_1,\lambda_2]^T$.  Using the minimum principle gives optimal open loop control
$$
u^*=argmin_{u\in R}H(X^*,u,\lambda)=-\eta_m^{-1}\lambda^{T}\left(
                                                                \begin{array}{c}
                                                                  0 \\
                                                                  1 \\
                                                                \end{array}
                                                              \right)
=-\eta_m^{-1}\lambda_2
$$
and optimal state equation
$$
\dot{X}^*=\nabla_{\lambda}H=\left(
          \begin{array}{c}
            y^* \\
            -(\alpha x^{*2}+\beta y^{*2}-\gamma)y^*-\omega^2x^*-\eta_m^{-1}\lambda_2\\
          \end{array}
        \right)
$$
with initial condition $X(t_k)=[x(t_k),\dot{x}(t_k)]^T$ and optimal costate equation
$$
\dot{\lambda}=-\nabla_{X}H=\left(
                               \begin{array}{c}
                                 \lambda_2(2\alpha x^*y^*+\omega^2) \\
                                 \lambda_2(\alpha x^{*2}+3\beta y^{*2}-\gamma)-\lambda_1-\theta_{\sigma}(y^*-r_{\sigma}) \\
                               \end{array}
                             \right)
$$
with the terminal condition
$$
\lambda(t_{k+1})=\left(
                   \begin{array}{c}
                     \theta_p(x^*(t_{k+1})-\hat{r}_p(t_{k+1})) \\
                     0 \\
                   \end{array}
                 \right)
$$
Let $\tilde{x}$ denote the approximation of the optimal solution $x^*$, then it is feasible to estimate the position error between the VP and the HP based on the collocation method as.
$$
|x^*-\hat{r}_p|=|x^*-\tilde{x}+\tilde{x}-\hat{r}_p|\leq |x^*-\tilde{x}|+|\tilde{x}-\hat{r}_p|
$$
Notice that $|x^*-\tilde{x}|$ is negligible due to the high approximation accuracy of numerical methods \cite{boor73}. In particular, considering that normally the optimal solution $x^*$ is not available, the approximate solution $\tilde{x}$ exactly corresponds to the position of the VP in the simulation. Thus, we mainly focus on the estimation of $|\tilde{x}-\hat{r}_p|$. For simplicity, we define $\tilde{x}(t)=a_0+a_1(t-t_k)+a_2(t-t_k)^2$, $\lambda_1(t)=b_0+b_1(t-t_k)+b_2(t-t_k)^2$ and $\lambda_2(t)=c_o+c_1(t-t_k)+c_2(t-t_k)^2$, where $a_i$, $b_i$ and $c_i$, $i\in\{0,1,2\}$ are unknown constants and $t\in [t_k,t_{k+1}]$. Substituting  $\tilde{x}(t)$, $\lambda_1(t)$ and $\lambda_2(t)$ into the above optimal state equation and costate equation at the boundary points yields the linear matrix equation
\begin{equation}\label{app_lin}
A_kX_k=B_k
\end{equation}
where
$$
A_k=\left(
  \begin{array}{ccccccccc}
    1 & 0 & 0 & 0 & 0 & 0 & 0 & 0 & 0 \\
    0 & 1 & 0 & 0 & 0 & 0 & 0 & 0 & 0 \\
    \theta_p & \theta_pT & \theta_pT^2 & -1 & -T & -T^2 & 0 & 0 & 0 \\
    0 & 0 & 0 & 0 & 0 & 0 & 1 & T & T^2 \\
    0 & 0 & 2 & 0 & 0 & 0 & \eta^{-1}_m & 0 & 0 \\
    0 & 0 & 0 & 0 & 1 & 0 & -(2\alpha x(t_k)y(t_k)+\omega^2) & 0 & 0 \\
    0 & 0 & 0 & 1 & 0 & 0 & -(\alpha x(t_k)^2+3\beta y(t_k)^2-\gamma) & 1 & 0 \\
    \theta_p & T\theta_p+\theta_{\sigma} & T(T\theta_p+2\theta_{\sigma}) & 0 & 0 & 0 & 0 & 1 & 2T \\
    0 & 0 & 0 & 0 & 1 & 2T & 0 & 0 & 0 \\
  \end{array}
\right)
$$
and
$$
B_k=\left(
                 \begin{array}{c}
                   x(t_k) \\
                   y(t_k) \\
                   \theta_p\hat{r}_p \\
                   0 \\
                   -(\alpha x(t_k)^2+\beta y(t_k)^2-\gamma)y(t_k)-\omega^2x(t_k) \\
                   0 \\
                   -\theta_{\sigma}(y(t_k)-r_{\sigma}(t_k)) \\
                   \theta_p\hat{r}_p+\theta_{\sigma}r_{\sigma}(t_{k+1}) \\
                   0 \\
                 \end{array}
               \right), \quad X_k=\left(
         \begin{array}{c}
           a_0 \\
           a_1 \\
           a_2 \\
           b_0 \\
           b_1 \\
           b_2 \\
           c_0 \\
           c_1 \\
           c_2 \\
         \end{array}
       \right)
$$
Solving equation (\ref{app_lin}) determines the vector of unknown constants
$$
X_k=A_k^{-1}B_k
$$
Thus, we obtain the approximate solution
\begin{equation*}
\tilde{x}(t)=x(t_k)+y(t_k)(t-t_k)+\frac{\mathcal{N}}{\mathcal{D}}(t-t_k)^2, \quad t\in [t_k,t_{k+1}]
\end{equation*}
where
\begin{equation*}
\begin{split}
\mathcal{N}=&2T[(\frac{r_{\sigma}(t_k)+r_{\sigma}(t_{k+1})}{2}-y(t_k))\theta_{\sigma}+(\hat{r}_p-x(t_k)-Ty(t_k))\theta_p]\\
&-\eta_m(\frac{T^2\omega^2}{2}+\alpha T^2 x(t_k)y(t_k)+\alpha Tx(t_k)^2+3\beta T y(t_k)^2-\gamma T+2)\\
&\cdot [(\alpha x(t_k)^2+\beta y(t_k)^2-\gamma)y(t_k)+\omega^2x(t_k)]
\end{split}
\end{equation*}
and
$$
\mathcal{D}=2T^2(\theta_pT+\theta_{\sigma})+2\eta_m(\frac{T^2\omega^2}{2}+\alpha T^2 x(t_k)y(t_k)+\alpha Tx(t_k)^2+3\beta T y(t_k)^2-\gamma T+2)
$$
Then we can compute
\begin{equation}\label{error}
\begin{split}
|\tilde{x}(t_{k+1})-\hat{r}_p(t_{k+1})|&=\lim_{t\rightarrow t_{k+1}}|\tilde{x}(t)-\hat{r}_p(t_{k+1})|\\
&=|x(t_k)+Ty(t_k)+\frac{\mathcal{N}}{\mathcal{D}}T^2-\hat{r}_p(t_{k+1})|\\
&\leq T^2(1-\theta_p)\frac{|2(x(t_k)-\hat{r}_p(t_{k+1}))+T(r_{\sigma}(t_k)+r_{\sigma}(t_{k+1}))|}{|\mathcal{D}|}+\eta_m\frac{|\mathcal{L}\cdot\mathcal{M}|}{|\mathcal{D}|}
\end{split}
\end{equation}
where
$$
\mathcal{L}=\frac{T^2\omega^2}{2}+\alpha T^2 x(t_k)y(t_k)+\alpha Tx(t_k)^2+3\beta T y(t_k)^2-\gamma T+2
$$
and
$$
\mathcal{M}=2(x(t_k)+Ty(t_k)-\hat{r}_p(t_{k+1}))-T^2(y(t_k)(\alpha x(t_k)^2+\beta y(t_k)^2-\gamma)+\omega^2x(t_k))
$$
Since $\hat{r}_p$, $r_{\sigma}$, $\mathcal{D}$, $\mathcal{L}$ and $\mathcal{M}$ are all bounded, it follows from inequality (\ref{error}) that the bound on the tracking error $|\tilde{x}(t_{k+1})-\hat{r}_p(t_{k+1})|$ converges to $0$ as $\theta_p\rightarrow 1$ and $\eta_m\rightarrow 0$.
Similarly, we can estimate the velocity error between the VP and the reference signal encoding the desired signature as follows
\begin{equation}\label{error_vel}
\begin{split}
|\dot{\tilde{x}}(t_{k+1})-r_{\sigma}(t_{k+1})|&=\lim_{t\rightarrow t_{k+1}}|\dot{\tilde{x}}(t)-r_{\sigma}(t)|\\
&=|y(t_k)+\frac{2\mathcal{N}}{\mathcal{D}}T-r_{\sigma}(t_{k+1})|\\
&\leq (1-\theta_{\sigma}) \frac{2T^2|T(y(t_k)-r_{\sigma}(t_{k+1}))+2(\hat{r}_p(t_{k+1})-x(t_k)-Ty(t_k))|}{|\mathcal{D}|}\\
&~~~~+\theta_{\sigma}\frac{2T^2|r_{\sigma}(t_k)-y(t_k)|}{|\mathcal{D}|}+2\eta_m\frac{|\mathcal{L}\cdot\mathcal{P}|}{|\mathcal{D}|}
\end{split}
\end{equation}
where
$$
\mathcal{P}=y(t_k)-r_{\sigma}(t_{k+1})-T[(\alpha x(t_k)^2+\beta y(t_k)^2-\gamma)y(t_k)+\omega^2x(t_k)]
$$
According to inequality (\ref{error_vel}), the bound of the velocity error goes to $0$ if $\theta_{\sigma}\rightarrow1$, $\eta_m\rightarrow 0$ and $r_{\sigma}(t_k)=y(t_k)$.

\end{document}